\documentclass[preprint,12pt,3p]{elsarticle}

\usepackage{booktabs}
\usepackage{lineno}
\modulolinenumbers[5]
\usepackage{tabularx}
\usepackage{graphicx}
\usepackage[ruled,noend,norelsize]{algorithm2e}
\usepackage{subfig}
\usepackage{amsmath}
\usepackage{amsthm}
\usepackage{amssymb}
\usepackage{bm}
\usepackage[table]{xcolor}
\usepackage{xspace} 
\usepackage{natbib} 
\usepackage{multirow}
\usepackage{setspace}
\usepackage{tikz}
\usepackage{float} 
\definecolor{yellow}{rgb}{1,1,0}
\definecolor{LightCyan}{rgb}{0.7,1,1}
\usepackage{pgf,tikz,pgfplots}
\pgfplotsset{compat=1.15}
\usepackage{mathrsfs}
\usetikzlibrary{arrows}
\usepackage[bookmarks=false,colorlinks]{hyperref}
\newtheorem{modl}{Model}
\newtheorem{prop}{Proposition}
\newtheorem{defi}{Definition}

\newenvironment{model}{\begin{samepage}\begin{modl}}{\end{modl}\end{samepage}}
\newenvironment{braced}{$\left\{\begin{minipage}{\columnwidth}}{\end{minipage}\right.$}

\journal{Elsevier}

\begin{document}

\begin{frontmatter}


\title{Freeway network design with exclusive lanes for automated vehicles under endogenous mobility demand}

\author[UNSWaddress]{Shantanu Chakraborty}
\author[UNSWaddress]{David Rey\corref{mycorrespondingauthor}}
\cortext[mycorrespondingauthor]{Corresponding author}
\ead{d.rey@unsw.edu.au}
\author[UMNaddress]{Michael W. Levin}
\author[UNSWaddress]{S. Travis Waller}

\address[UNSWaddress]{School of Civil and Environmental Engineering, UNSW Sydney, Sydney, NSW, 2052, Australia}
\address[UMNaddress]{Department of Civil, Environmental, and Geo-Engineering, University of Minnesota, Minneapolis, MN, 55455, United States}

\begin{abstract}
Automated vehicles (AV) have the potential to provide cost-effective mobility options along with overall system-level benefits in terms of congestion and vehicular emissions. Additional resource allocation at the network level, such as AV-exclusive lanes, can further foster the usage of AVs rendering this mode of travel more attractive than legacy vehicles (LV). However, it is necessary to find the crucial locations in the network where providing these dedicated lanes would reap the maximum benefits. In this study, we propose an integrated mixed-integer programming framework for optimal AV-exclusive lane design on freeway networks which accounts for commuters' demand split among AVs and LVs via a logit model incorporating class-based utilities. We incorporate the link transmission model (LTM) as the underlying traffic flow model due to its computational efficiency for system optimum dynamic traffic assignment. The LTM is modified to integrate two vehicle classes namely, LVs and AVs with a lane-based approach. The presence of binary variables to represent lane design and the logit model for endogenous demand estimation results in a nonconvex mixed-integer nonlinear program (MINLP) formulation. We propose a Benders' decomposition approach to tackle this challenging optimization problem. Our approach iteratively explores possible lane designs in the Benders' master problem and, at each iteration, solves a sequence of system-optimum dynamic traffic assignment (SODTA) problems which is shown to converge to fixed-points representative of logit-compatible demand splits. Further, we prove that the proposed solution method converges to a local optima of the nonconvex problem and identify under which conditions this local optima is a global solution. The proposed approach is implemented on three hypothetical freeway networks with single and multiple origins and destinations. \textcolor{blue}{Our numerical results reveal that the optimal lane design of freeway network is non-trivial while accounting for endogenous demand of each mode.}
\end{abstract}


\end{frontmatter}

\renewcommand\thesection{\arabic{section}}
\newcommand{\x}{x_{i,t}}
\newcommand{\y}{y_{i,j,lv}^{k}(t)}
\newcommand{\yav}{y_{i,j,av}^{k}(t)}
\newcommand{\zp}{z^{k+}_{i,lv}(t)}
\newcommand{\zpav}{z^{k+}_{i,av}(t)}
\newcommand{\zn}{z^{k-}_{i,lv}(t)}
\newcommand{\znlv}{z^{k-}_{i,lv}(t)}
\newcommand{\znav}{z^{k-}_{i,av}(t)}
\newcommand{\lvd}{d_{i,lv}^{k,t}}
\newcommand{\avd}{d_{i,av}^{k,t}}
\newcommand{\origd}{d_{i}^{0,k,t}}
\newcommand{\klv}{K_{jam}}
\newcommand{\kav}{K_{jam}}
\newcommand{\bigzp}{Z_{i}^{k,t+}}
\newcommand{\bigzn}{Z_{i}^{k,t-}}
\newcommand{\demlv}{D_{lv}^{o,d}(t)}
\newcommand{\demav}{D_{av}^{o,d}(t)}
\newcommand{\problv}{p_{lv}^{o,d}}
\newcommand{\probav}{p^{o,d}}
\newcommand{\totdem}{D^{o,d}(t)}
\newcommand{\bxij}{\bar{x}_{ij}}
\newcommand{\xij}{x_{ij}}
\newcommand{\xji}{x_{ji}}
\newcommand{\tij}{t_{ij}}
\newcommand{\tji}{t_{ji}}
\newcommand{\lij}{l_{ij}}
\newcommand{\lji}{l_{ji}}
\newcommand{\ai}{\alpha_i}
\newcommand{\bi}{\beta_i}
\newcommand{\aj}{\alpha_j}
\newcommand{\bj}{\beta_j}

\newcommand{\Ydom}{\mathcal{Y}}
\newcommand{\Zdom}{\mathcal{Z}}
\newcommand{\Bdom}{\mathcal{B}}
\newcommand{\Tdom}{\mathcal{T}}
\newcommand{\Pdom}{\mathcal{P}}
\newcommand{\PF}{\texttt{FNDP}}
\newcommand{\MP}{\texttt{MP}}
\newcommand{\SP}{\texttt{SP}(\bm{b},\bm{p})}
\newcommand{\ulv}{U_{lv}^{o,d}}
\newcommand{\uav}{U_{av}^{o,d}}
\newcommand{\tlv}{\tau_{lv}^{o,d}}
\newcommand{\tav}{\tau_{av}^{o,d}}
\newcommand{\blv}{\beta^{o,d}_{\tau_{lv}}}
\newcommand{\bav}{\beta^{o,d}_{\tau_{av}}}

\SetArgSty{textnormal}

\section{Introduction}
As mobility demand increases, transport infrastructures remain under-utilized majorly due to various human factors e.g., slow reflexes, poor and heterogeneous driving behaviour, safety concerns etc. Automation in transport sector is gaining more attention than ever to minimize these human inputs to exploit transport resources in an efficient and sustainable way. This growing attention in automation translates to a projected growth of US\$173B in global autonomous driving market by 2030 \cite{frost&sullivan}. Potential benefits of automated vehicles (AV) include improved throughput \cite{shladover2012impacts,vander2017operational,levin2016cell}, traffic safety \cite{fernandes2012platooning}, travel speed, energy consumption \cite{mersky2016fuel} and vehicular emissions \cite{greenblatt2015autonomous}. 

Regardless of the benefits of AVs, one would be too naive to assume that AVs will be immediately adopted by legacy vehicle (LV) owners in near future. A more reasonable assumption would be the existence of a transition period where interactions between LVs and AVs exist, leading to a gradual increment in market penetration of AVs. This transition period will be crucial as safety might be compromised in mixed operations of LVs and AVs, especially in case of arterial networks, involving pedestrians, cyclists and signalized intersections where fewer AVs have been found to have a negative impact on the average travel speed and string stability of the traffic flow \cite{van2006impact,talebpour2016influence,monteil2018mathcal}.

\textcolor{blue}
{These negative impacts might have resulted due to the behavioural response of the drivers of LVs unfamiliar with vehicles with automation. To minimize these negative impacts, it is necessary to introduce AVs gradually in a network starting with the less conflicting regions such as freeways, hence allowing the technology to mature and users to adopt to changing driving conditions induced by the presence of automated vehicles. During this transition period, safety concerns due to mixed-vehicular interactions could be further reduced via AV-exclusive lanes which separate AVs and LVs in a network, providing a smooth transition of the current transport system to an automation-heavy transport system. These dedicated lanes would be advantageous in the current situation of imperfect AV technologies, which pose critical problems in traffic safety during lane changing in a congested network with mixed traffic.}

\subsection{\textcolor{blue}{AV-exclusive lane design on freeways}}
\citet{yu2019impact} adopted a microscopic traffic simulation method to investigate the efficiency and safety of mixed traffic on highways with AV-exclusive lanes. Although, the safety of mixed traffic was found to be worsened with low market penetration of AVs, the simulation results from the car-following models of this study showed an increment of up to 84\% in throughput of the traffic network due to presence of AV-exclusive lanes. These dedicated lanes may also facilitate cooperative adaptive cruise control (CACC) contributing to better traffic flow performance \cite{van2006impact}, improved highway capacity \cite{milanes2013cooperative}, decrease in fuel consumption and emissions \cite{ploeg2011design}. On freeways, CACC was found to significantly increase capacity with a moderate to high market penetration rate \cite{shladover2012impacts}. {However, with a naive deployment strategy of AV-exclusive lanes, CACC was found to increase the total system travel time on freeways \cite{MELSON2018114}}. \citet{alireza2017} simulated traffic flow under different penetrations rates of AVs on a two-lane and a four-lane freeway segment in Chicago, Illinois with mandatory and optional usage of AV-exclusive lanes. At market penetration rates of more than 50\% for the two-lane highway and 30\% for the four-lane highway, a potential benefit in terms of throughput and travel time reliability was observed in this study, with optional use of the AV-exclusive lanes yielding the most benefit. In another study, this optional use of AV-exclusive lanes was found to be beneficial only with a high market share of AVs \cite{mahmassani201650th}. 

\subsection{{AV market share}}
With the advent of AV technology in the market, it is critical to gain individual motivations for choosing to own or use AVs as a service. \citet{hyland2018dynamic} modelled a sequentially shared AV-service as a dynamic AV-traveller assignment problem where advanced optimization-based assignment strategies were found to reduce the average traveller wait times with a smaller fleet size relative to the demand rate. Based on a stated preference survey from 721 individuals living across Israel and North America, 44\% of the sample population were found to prefer LVs over AVs \cite{haboucha2017user}. Another stated choice survey for the adoption of shared AVs suggested service attributes including travel cost, travel time and waiting time to be critical determinants of the use of shared AVs \cite{krueger2016preferences}. Hence, capturing the effects of the perceived benefits of AVs on AV-adoption or using AVs as a service becomes imperative while modelling the multi-modal traffic flow with LVs and AVs. A few studies incorporated demand models with AV-exclusive lanes. \citet{chen2016optimal} proposed a model to optimally deploy AV lanes considering the endogenous AV market penetration and a multi-class network equilibrium model. Here, a diffusion model was developed to forecast the evolution of AV market penetration over time. Another study by \citet{movaghar2020optimum} modelled optimum location of AV-exclusive lanes with link capacity being a function of AV proportion in that link. Several other studies on AV-exclusive lanes with AV-penetration rate found these lanes to improve overall throughput \cite{ye2018impact,ma2019influence} on general and freeway networks. At network intersections, reservation-based models have been shown to have the potential to increase intersection capacity \cite{dresner2004multiagent,fajardo2011automated,qian2014priority,levin2017conflict}. More recently, \citet{rey2019blue} proposed a hybrid network control policy in such networks with dedicated lanes that provide access to ``blue phases" during which only AVs can traverse traffic intersections.     

To allow the transportation system to adapt with such disruptive technologies, we explore the potential benefits of introducing AV-exclusive lanes on freeways where the conflict points are considerably lesser than in local arterial networks and require minimal technological interventions. From previous studies, it is apparent that deploying AV-exclusive lanes in all links in a network may not be advantageous as it may significantly increase congestion and reduce throughput in the regular lanes; yet no systematic approach to identify the optimal location of AV-exclusive lanes has been proposed. In this work, we propose to investigate the problem of locating AV-exclusive lanes on freeway networks so as to reap maximum benefits of the deployed infrastructure. 

\subsection{{System Optimum Dynamic Traffic Assignment (SODTA)}}
We build on the literature of system optimum dynamic traffic assignment (SODTA) formulations to represent traffic dynamics, and obtain an analytical model amenable for exact optimization. Although a few studies \cite{long2019link,long2018link,levin2017congestion,shantanu2018} proposed SODTA approaches to model user route choice in general networks, this behavioral assumption is limiting, especially in the context of endogenous travel demand. Alternative approaches such as cross-nested and user-equilibrium models have also been proposed for characterizing different route choice principles \cite{wang2019multiclass,MELSON2018114}, these research efforts are focused on static traffic assignment formulations or do not consider endogenous travel demand. To circumvent such assumptions, we focus on the design of freeway networks with single path per OD pair and we leave the study of general networks for future research. In this context, there are no route choice decisions, and the SODTA component of the formulation provides an analytical representation of traffic dynamics in the freeway network which maintains first-in-first-out (FIFO) conditions. Further, we do not model lane-changing behaviour of vehicles and do not account for vehicle holding in the network which is defined as the reluctance of vehicles moving forward from upstream to downstream links even with an availability of vacant spaces downstream \cite{long2018link,long2019link}. 

For modelling the traffic flow in the network, we adopt the link transmission model (LTM) \cite{yperman2005link} and modify it to account for two vehicle classes: LVs and AVs. The LTM has been found to be useful to scale-up SODTA formulations for network design \cite{shantanu2018} as well as for solving routing problem for shared AVs in congested networks \cite{levin2017congestion}. In the proposed formulation, we avoid mixed vehicular interactions of LVs and AVs by providing AV-exclusive lanes allowing AVs to make use of their automated features only on AV-exclusive lanes and assume both LVs and AVs to behave identically on regular lanes. These automated features include level 1 automation on the Society of Automotive Engineers automation scale \citep{sae2013definitions} with cooperative adaptive cruise-control (CACC), speed harmonisation and cooperative merging. We capture the demand split of each vehicle class with a logit model, embedded in the proposed formulation. We model lane design decisions using binary variables, thus the proposed formulation results in a challenging non-convex MINLP. We adopt a Benders' decomposition approach to implement this MINLP on multi-OD freeway networks which embeds a fixed-point algorithm in its sub-problem to account for endogenous travel demand.

\subsection{{Our contributions}}
{
This study makes the following contributions.
We propose a formulation to identify the optimal design of dedicated lanes for AVs on a freeway network wherein the demand of LVs and AVs is endogenous. To the best of our knowledge, this is the first integrated modelling framework that combines dynamic network loading with an endogenous choice model for AV-exclusive lane design problems.
Given the critical speculation of vehicle owners on AV adoption, a prolonged transition period is evident in the near future during which both LVs and AVs will coexist in road networks. To model such hybrid scenario, this study adapts the traditional LTM of traffic flow to account for multiple vehicle classes on a network with AV-exclusive lanes.
This multi-class traffic flow model along with the endogenous choice model and lane design variables transforms the traditional SODTA model with a linear structure into a non-convex and mixed-integer non-linear program. We propose a Benders' decomposition approach in which the sub-problem is formulated as a fixed-point problem and is solved by a sequence of linear programs under varying demand rates. We prove that the proposed freeway network design problem always admits at least one feasible solution and present an algorithm that is guaranteed to converge to a fixed-point solution.
We prove the existence of a fixed-point and our proposed algorithm always converges to this fixed-point. We also identify the conditions under which our algorithm is globally optimal. We implement the proposed algorithm on freeway networks which reveals that deploying a maximum of AV-exclusive lanes does not necessarily minimize network travel time.} 

\textcolor{blue}{The gradual advent of automation technologies in the transport industry is expected to provide enhanced coordination and cooperation between (network infrastructure) supply and (travel) demand in the near future. Hence, it is imperative to prepare the current transport system steadily to accept such disruptive technologies and reap the maximum benefits from them. The main purpose of this study is to contribute to this expedition towards future mobility. We advocate adopting a gradual implementation approach allowing the use of automation technologies in areas with less number of conflict points (e.g., freeways). To further dampen the effects of technological disruptions, this study explores the concept of AV-exclusive lanes in a comprehensive manner where vehicles with automation features are allowed to wield their technological advancements only on these exclusive lanes. These exclusive lanes are expected to shelve mixed-vehicular interactions between LVs and AVs allowing the technology to mature. As technology and user adoption evolve, mixed-vehicular interactions within lanes may emerge on additional elements of transportation networks.}

This paper is organized into five sections. Section \ref{problem} presents the problem formulation followed by the solution methodology in Section \ref{solution}. The proposed formulation is studied on two numerical networks, presented in Section \ref{numerical}. Section \ref{conclusion} presents the key findings of the study along with future research directions.  

\section{Freeway network design problem}
\label{problem}
In this study, we develop an MINLP model for optimal AV-exclusive lane design along with endogenous estimation of AV demand. We develop this MINLP in three stages. To begin with, we propose a linear programming (LP) framework for a lane-based LTM formulation with fixed AV-exclusive lanes, fixed proportion of AVs in the network and a system level objective. Further, we introduce a binary variable to obtain optimal lane design for AVs for an improved system performance, resulting in a mixed-integer linear program (MILP). Finally, we bring in a logit model to estimate the endogenous demand for each vehicle class which introduces non-linearity in the model, resulting in an MINLP. We circumvent this non-linearity by proposing a fixed-point algorithm along with method of successive averages (MSA) to obtain convergence of the fixed-point as explained in Section \ref{solution}. The logit model embedded in this fixed-point algorithm, is solved iteratively along with the multi-OD model of SODTA to find a proportion of AVs satisfying a mode-choice equilibrium. The entire formulation involving AV-exclusive lane allocation, the multi-OD SODTA and the fixed-point algorithm is decomposed with Benders' decomposition method to disentangle the binary lane allocation variables with the rest of the formulation. To the best of our knowledge, such an integrated framework has not been proposed in the literature. The subsequent sections present this framework in detail.

\subsection{{Problem definition}}
{We consider two vehicle classes: LVs and AVs. The LVs are regular vehicles without any connectivity features for vehicle-to-vehicle communication, whereas, AVs are equipped with connected and automated features such as CACC, speed harmonisation and cooperative merging}. These features correspond to level 1 automation on the Society of Automotive Engineers automation scale \cite{sae2013definitions}. In a mixed traffic network with these two vehicle classes, an AV can fully utilise its automated features only if the leading/following vehicles are AVs as well. To allow AVs to get full benefits of their automated features, this study proposes a lane-based approach with AV-exclusive lanes which forbid the entry of LVs. The regular lanes allow both the vehicle classes where AVs behave like LVs with restrained automated features. As providing more number of these exclusive lanes might affect the traffic flow of LVs, they are provided only at those crucial locations which would improve the overall system performance. The allocation of these exclusive lanes might impact the total demand in the network as well. However, in this study this total demand is assumed to be fixed with varying class-wise demand. \textcolor{blue}{The problem setting is a multi-OD freeway network with multiple paths between each OD pair where the route choice modelling is governed by the objective function of the SODTA.} Further, lane-changing behaviour and vehicle holding are not modelled in the proposed approach and are kept for future research.

{In this context, we define the freeway network design problem (FNDP) as follows.
\begin{defi}[FNDP]
	Given a freeway network with multiple origins and destinations and two vehicle classes, i.e. AV and LV, find the optimal location of AV-exclusive lanes on the network. 
\end{defi}
}

{We next introduce the mathematical formulation to model the proposed FNDP.} 

\subsection{{Model formulation}}
{In the FNDP, we categorise the lanes in the network into two groups: regular lanes and candidate AV lanes. The regular lanes do not participate in the lane allocation problem and provides at least one path for the LVs to reach their destinations. On the other hand, the candidate AV lanes could be converted to AV-exclusive lanes in the network through a binary decision variable if it improves the overall objective function that minimises the total system travel time (TSTT) involving both vehicle classes.}

Let $G = (N,A)$ be a directed network where $N$ is the set of nodes and $A$ is the set of arcs. In the proposed formulation, each node of the network is represented by a set of incoming and outgoing arcs. The set of arcs is partitioned into three subsets: source centroid connectors denoted as $A_r$, sink centroid connectors as $A_s$, and physical links formed by the remaining links of the set; i.e., $A \setminus \{A_r \cup A_s\}$. The set of origin-destination pairs is represented by $W$. $\Gamma^-(i)$ and $\Gamma^+(i)$ represent the set of predecessor and successor links of link $i$. Table \ref{tab:notations} presents the rest of the notations of the proposed formulation.
\begin{table}[!h]
	\caption{Notations}
	\resizebox{\textwidth}{!}{
		\begin{tabular}{@{}ll@{}}
			\toprule
			\textit{Sets}                                                  &                                                                                                                                                           \\ \midrule
			$A$                                                                             & set of all lanes and centroid connectors                                                                                                                  \\
			$A_r$                                                                           & set of source centroid connectors                                                                                                                         \\
			$A_s$                                                                           & set of sink centroid connectors                                                                                                                           \\
			$A_c$                                                                           & set of source centroid connectors and physical lanes                                                                                                                          \\
			$A_{av}$                                                                        & set of candidate AV-exclusive lanes                                                                                                                           \\
			$M$ 																			& set of links
			\\
			$W$                                                                             & set of origin-destination pairs                                                                                                                           \\
			$\Gamma^-(i)$                                                                   & set of predecessor lanes of lane $i \in A$                                                                                                                      \\
			$\Gamma^+(i)$                                                                   & set of successor lanes of lane $i \in A$                                                                                                                        \\
			$T$                                                                             & set of discretized time steps for traffic flow propagation ($t_0,t_1,...,t_n$)                                                                            \\
			$R$                                                                            & set of discretized time steps for demand loading\\
			$\Pi^{o,d}$ & set of paths between OD pair $(o,d)$\\
			$m(l)$ & set of lanes in link $l$ \\                                                                                                                                                                                                            \midrule                                                                                                          
			\begin{tabular}[c]{@{}l@{}}\textit{Parameters}\end{tabular} & 
			\\ \midrule                                                                                                                                                          
			$\totdem$                                                                       & total demand from $o$ to $d$ at $t$\\
			$q_{lv}$                                                                        & capacity of a regular lane                                                                                                                                    \\
			$q_{av}$                                                                        & capacity of an AV-exclusive lane                                                                                                                                    \\
			$L_i$                                                                           & length of lane $i \in A$                                                                                                                                          \\
			$\klv$                                                                          & jam density                                                                                                                                               \\
			$v_{f}$                                                                         & free-flow speed                                                                                                                                           \\
			$w_{lv}$                                                                     & backward shockwave speed on a regular lane                                                                                                                    \\
			$w_{av}$                                                                        & backward shockwave speed on an AV-exclusive lane                                                                                                                    \\
			$\delta$                                                                        & discretized time step for traffic flow propagation                                                   			\\
			$\beta$                                                                   & total amount of utility gained while making a trip                                                                                                         \\
			$\blv$                                                                   & disutilities per unit travel time of LVs for $(o,d) \in W$                                                                                                         \\
			$\bav$                                                                   & disutilities per unit travel time of AVs for $(o,d) \in W$ \\                                                                                                                                                       \midrule                                                                                                        
			\textit{Variables}\\ 
			\midrule                                            
			$\y \geq 0$                                                                            & transfer flow of LVs from lane $i \in A$ to lane $j \in A$ destined to $k \in A_s$ at time $t \in T$                                                                                \\
			$\yav \geq 0$                                                                          & transfer flow of AVs from lane $i \in A$ to lane $j \in A$ destined to $k \in A_s$ at time $t \in T$                                                                                \\
			$\zp \geq 0$                                                                           & cumulative inflow of LVs on lane $i \in A \setminus A_s$ destined to $k \in A_s$ at time $t \in T$                                                                        \\
			$\zpav \geq 0$                                                                         & cumulative inflow of AVs on lane $i \in A \setminus A_s$ destined to $k \in A_s$ at time $t \in T$                                                                       \\
			$\zn \geq 0$                                                                           & cumulative outflow of LVs on lane $i \in A \setminus A_s$ destined to $k \in A_s$ at time $t \in T$                                                                       \\
			$\znav \geq 0$                                                                         & cumulative outflow of AVs on lane $i \in A \setminus A_s$ destined to $k \in A_s$ at time $t \in T$                                                                      \\
			$b_i \in \{0,1\}$                                                                           & binary variable indicating whether a lane $i \in A_{av}$ is AV-exclusive (1) or not (0)                                                                                                                        \\
			$\probav \in [0,1]$                                                                       & fraction of AV demand for $(o,d) \in W$                                                                              \\
			$\tau_{lv,l} \geq 0$                                                                   & average link travel time for LVs on link $l$                                                                                                         \\
			$\tau_{lv,l} \geq 0$                                                                   & average link travel time for AVs on link $l$
			\\
			$\tlv \geq 0$                                                                   & average path travel time for LVs for $(o,d) \in W$                                                                                                         \\
			$\tav \geq 0$                                                                   & average path travel time for AVs for $(o,d) \in W$
			\\ \bottomrule
	\end{tabular}}
	\label{tab:notations}
\end{table}

\subsubsection{{Endogenous demand model}}
\label{demand}
To estimate the AV proportion in the network based on the gain/loss in utilities of AVs compared to the LV mode, a logit model is adopted for mode choice. The total vehicular demand in the network is the sum of LVs and AVs in the network as follows.

We denote the time-varying total vehicle demand at time $t$ between each origin-destination pair in the network by $\totdem$ and assume it to be fixed in our formulation. $\totdem$ is presented as a sum of the demands of two vehicle classes, LVs and AVs, in Eq. \eqref{totdem}.
\begin{align}
\label{totdem}
\totdem &= \demlv+\demav && \forall (o,d) \in W, \forall t \in R
\end{align}

Though, the total demand in the network is fixed, the demand corresponding to each vehicle class varies depending on the proportions of AVs ($\probav$) between each OD pair in the network and they are obtained from Eqs. \eqref{demlv} and \eqref{demav}. 
\begin{subequations}
	\begin{align}
	\label{demlv}
	\demlv&=\totdem(1-\probav) && \forall (o,d) \in W, \forall t \in R\\
	\label{demav}
	\demav&=\totdem\probav && \forall (o,d) \in W, \forall t \in R
	\end{align}
\end{subequations}

In the proposed formulation, the demand for each vehicle class is endogenous to the proportions of AVs ($\probav$), which is obtained based on the attractiveness of the modes in the network. We adopt a logit model to quantify this attractiveness, as summarized in Eq. \eqref{probav}.
\begin{align}
\label{probav}
\probav = \frac{e^{\uav}}{e^{\ulv}+e^{\uav}}
\end{align}

We assume that the utility of each mode depends only on the average travel times of all the vehicles of that mode between each OD pair. The average travel times of LV and AV are denoted by $\tlv$ and $\tav$ and the utility of each mode is obtained from Eqs. \eqref{ulv} and \eqref{uav}.
\begin{subequations}
	\begin{align}
	\label{ulv}
	\ulv &= \beta - \blv\tlv &&\forall (o,d) \in W\\
	\label{uav}
	\uav &= \beta - \bav\tav &&\forall (o,d) \in W
	\end{align}
\end{subequations}

Here, $\beta$ represents the total amount of utility gained while making a trip, whereas, $\blv$ and $\bav$ are the disutilities created per unit travel time by LVs and AVs respectively. We refer to a study by \citet{neeraj2018} to obtain the values of $\blv$ and $\bav$.

We modify the logit model presented in Eq. \eqref{probav} with the difference in utilities between the modes where the utility functions consist of the coefficients ($\blv,\bav$) of average travel times of corresponding modes, as presented in Eq. \eqref{pav}.
\begin{align}
\label{pav}
\probav &= \frac{1}{e^{\big(\beta^{o,d}_{\tau_{av}}\tau^{o,d}_{av}-\beta^{o,d}_{\tau_{lv}}\tau^{o,d}_{lv}\big)}+1} &&\forall (o,d) \in W
\end{align}

The OD average link travel times of each vehicle class are obtained from the lane-based formulation of LTM which is the underlying traffic flow model in the proposed formulation. {Here, we assume that the proportion of AVs ($\probav$) and average travel times ($\tlv, \tav$) are not time-dependent as a new $\probav$ and lane-allocation strategy at each timestep may not be realistic.}

The mode choice model presented in this study involves simplistic utility functions with mode-wise average freeway network travel times. We acknowledge that the commuters' decision to choose AV mode for travel might depend on other factors such as cost of AV adoption as well as travel cost, including maintenance and energy, which are not considered in this study while estimating the utilities. \textcolor{blue}
{However, freeway network travel times could be a major factor in deciding the overall AV proportion for an OD pair. If buying an AV allows travellers to save considerable travel time on their daily commute, it might be considered a strong incentive.}

\textcolor{blue}{On the other hand, the proposed formulation is flexible enough to incorporate a more refined nested logit, or probit model. Since, there exists a dedicated inner-loop to solve the choice model, the content of that loop could be replaced with advanced choice models that either truly capture AV ownership and/or choice models that are able to capture the multiple facets of automated mobility, e.g. AV as a service.} We emphasise that incorporating such refinements within the utility model can be achieved without altering the mathematical structure of the proposed optimisation approach, hence we leave this for future research.

\subsubsection{{Dynamic network model}}
\label{ltm}

The LTM, proposed by \citet{yperman2005link}, is a numerical solution method for dynamic network loading, developed based on the first order kinematic wave theory \cite{lighthill1955kinematic,richards1956shock}. In this study, LTM is chosen as it involves fewer variables per link compared to models such as cell transmission model (CTM) \cite{daganzo1994cell}. We adapt the conventional LTM to accommodate two vehicle classes by introducing two types of lanes: AV-exclusive and regular lanes. An AV-exclusive lane differs from a regular lane in terms of following headway, capacity and speed of backward shockwave propagation. The difference in these traffic flow characteristics affects the fundamental diagrams of traffic flow significantly as explained in the following subsection. 

\subsubsection*{2.2.2.1 Fundamental diagram}
The fundamental diagram of traffic flow reflects the relationship among the macroscopic traffic flow parameters of a network: traffic flow, density and speed. These relationships approximate all possible stationary traffic states during the analysis period and provide significant insights regarding the overall behaviour of traffic in a network.

The shape of the fundamental diagrams depicting these relationships may vary depending on the assumptions and approximations of a study. \citet{greenshields} was the first to propose a parabolic relationship between traffic flow and density. Later on, \citet{newell1993simplified} provided a simplified approach to the kinematic wave theory of traffic flow and developed a triangular shaped fundamental diagram, defined by three parameters: a fixed free-flow speed ($v_f$), capacity or maximum flow ($q$) and a jam density ($K_{jam}$). \citet{yperman2005link} adopted this simplified fundamental diagram while developing the LTM which is the underlying traffic flow model in our formulation. 

In a network with AV-exclusive lanes, the macroscopic traffic flow parameters may be significantly affected by faster reaction times of AVs leading to reduced following headway, increased throughput and faster propagation of backward shockwave due to congestion. \citet{LEVIN2016103} found considerable difference in the shape of the fundamental diagram for different reaction times of a characteristic vehicle. This study also showed how capacity and wave speed increase as the AV proportion increases with the human drivers having double the reaction time of AVs. \citet{tientrakool2011highway} demonstrated that due to tighter time and space headways among vehicles, the capacity of a lane could be approximately tripled by converting it into an AV-exclusive lane. Hence, while comparing traffic flow on AV-exclusive and regular lane, the shape of the triangular fundamental diagram will be significantly different due to the changes in capacity ($q$) and backward wave speed ($w$) leading to the same jam density ($K_{jam}$).

\begin{figure}[H]
	\centering
	\includegraphics[width=0.7\textwidth]{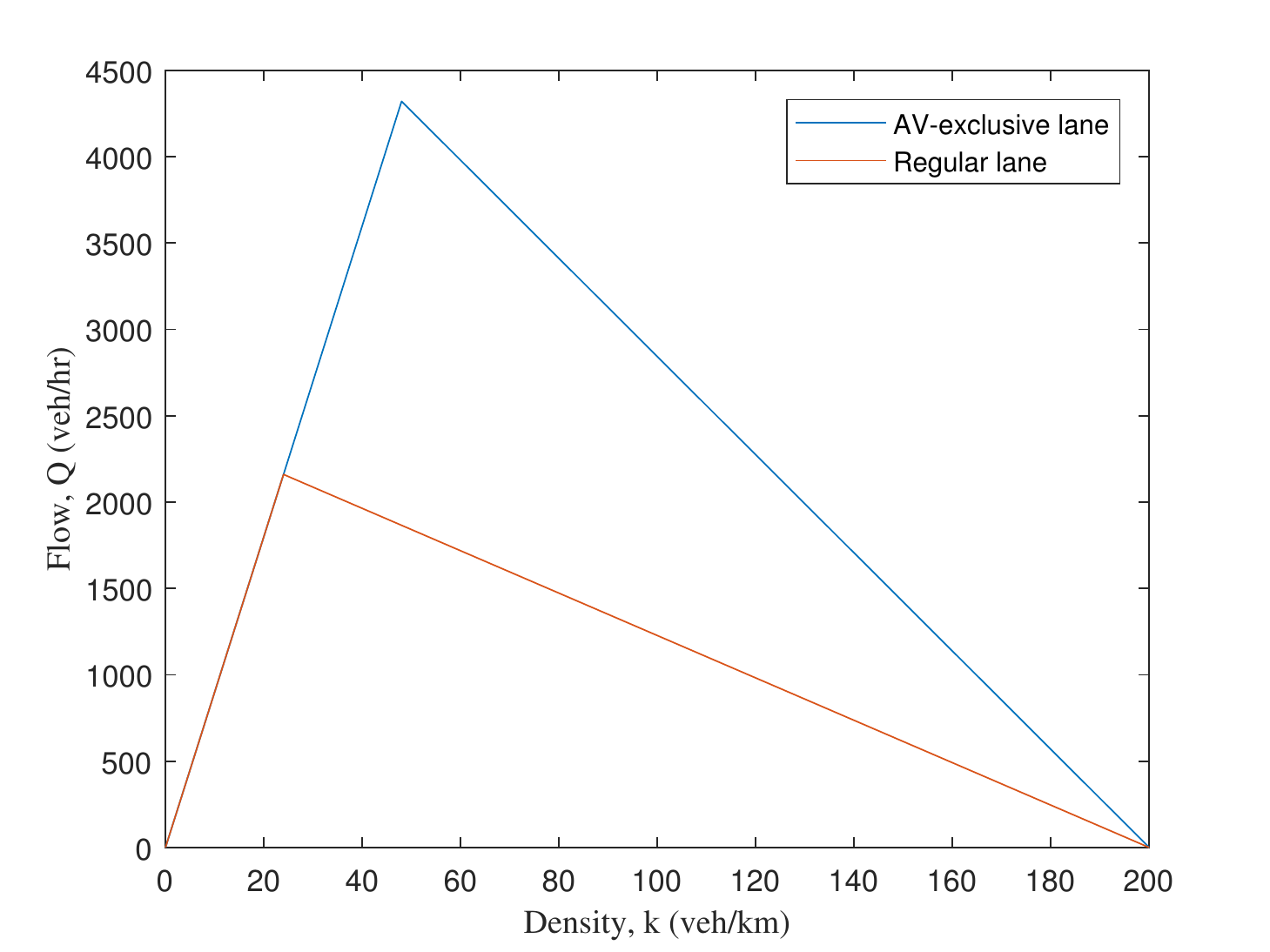}
	\caption{Fundamental diagrams of traffic flow for two lane types}
	\label{fig:fd}
\end{figure}

In this study, we consider the AV-exclusive lanes to have double the capacity ($q_{av}=2q_{lv}$) and backward wave speed ($w_{av}=2w_{lv}$) of that of the regular lane while keeping the free-flow speed ($v_f$) and jam density ($K_{jam}$) the same for both lane types. The fundamental diagrams of traffic flow on both of these lane types is shown in Figure \ref{fig:fd}.

\subsubsection*{2.2.2.2 Traffic flow propagation}
The LTM keeps track of the vehicular flow in the network with cumulative inflows and outflows of each link at each time-step. We develop a lane-based LTM where the vehicle class-specific cumulative inflows and outflows from each lane $i$ towards destination $k$ at time $t$ are denoted by $\zp \big(\zpav\big)$ and $\zn \big(\znav\big)$ respectively. The demand corresponding to each vehicle class is loaded into the network as the cumulative inflow to the source centroid connectors as shown in Eqs. \eqref{demandlv} and \eqref{demandav}.

\begin{subequations}
	\begin{align}
	\label{demandlv}
	\zp &= \sum_{t'<t}D^{i,k}(t')(1-p^{i,k}) \qquad &&\forall i \in A_r, \forall (i,k) \in W, \forall t \in T\\
	\label{demandav}
	\zpav &= \sum_{t'<t}D^{i,k}(t')p^{i,k} \qquad &&\forall i \in A_r, \forall (i,k) \in W, \forall t \in T
	\end{align}
\end{subequations}

The cumulative inflow to the other lanes at time $t$ is defined as the sum of transfer flows from all the incoming lanes predecessor ($\Gamma^-$) to that lane over all the timesteps up until $t$. We define these in Eqs. \eqref{zplv} and \eqref{zpav}.
\begin{subequations}
	\begin{align}
	\label{zplv}
	\zp &= \sum_{t'< t} \sum_{h \in \Gamma^-(i)}y_{h,i,lv}^{k}(t') \qquad &&\forall i \in A\setminus \{A_r,A_s\}, \forall k \in A_s, \forall t \in T\\
	\label{zpav}
	\zpav &= \sum_{t'< t} \sum_{h \in \Gamma^-(i)}y_{h,i,av}^{k}(t') \qquad &&\forall i \in A\setminus \{A_r,A_s\}, \forall k \in A_s, \forall t \in T
	\end{align}	
\end{subequations}

Similarly, the cumulative outflows ($\zn,\znav$) from a lane at time $t$ is defined as the sum of transfer flows to all the outgoing lanes successor ($\Gamma^+$) to that lane over all the timesteps up until $t$. We define these in Eqs. \eqref{znlv} and \eqref{znav}.
\begin{subequations}
	\begin{align}
	\label{znlv}
	\zn &= \sum_{t'< t} \sum_{j \in \Gamma^+(i)}y_{i,j,lv}^{k}(t') \qquad &&\forall i \in A\setminus A_r, \forall k \in A_s, \forall t \in T\\
	\label{znav}
	\znav &= \sum_{t'< t} \sum_{j \in \Gamma^+(i)}y_{i,j,av}^{k}(t') \qquad &&\forall i \in A\setminus A_r, \forall k \in A_s, \forall t \in T
	\end{align}
\end{subequations}

The LTM has been built based on three flow components: sending flow, receiving flow and transfer flow. Sending flow is defined as the amount of vehicular flow allowed to go out from link $i$ to link $j$ respecting its flow capacity. \citet{yperman2005link} derived the equation of sending flow based on the propagation of a free-flow traffic state at the upstream boundary of a link transmitting to the downstream boundary $\frac{L_i}{v_{f,i}}$ (link free-flow travel time) time units later. In our lane-based LTM formulation, we implement this concept for each of the vehicle classes as presented in Eqs. \eqref{sendlv} and \eqref{sendav}.

\begin{subequations}
	{
	\begin{multline}
	\label{sendlv}
	\sum_{j \in \Gamma^+(i)}\y \leq \big(z_{i,lv}^{k+}(t_s) - \znlv\big) \qquad \forall i \in A\setminus A_s, \forall k \in A_s, \forall t \in T\setminus \{t_n\} \\\text{ where, }t_s = t+\delta-\frac{L_i}{v_{f}}
	\end{multline}
	}
	\begin{multline}
	\label{sendav}
	\sum_{j \in \Gamma^+(i)}\yav \leq \big(z_{i,av}^{k+}(t_s) - \znav\big) \qquad \forall i \in A\setminus A_s, \forall k \in A_s, \forall t \in T\setminus \{t_n\} \\\text{ where, }t_s = t+\delta-\frac{L_i}{v_{f}}
	\end{multline}
\end{subequations}

Eqs. \eqref{send2lv} and \eqref{send3lv} present the capacity constraints on sending flow for regular and candidate AV lanes respectively. On regular lanes, the total flow of LVs and AVs is restricted to the capacity of regular lanes ($\delta q_{lv}$) whereas, on candidate AV lanes, we introduce a binary variable $b_i \in \{0,1\}$, $\forall i \in A_{av}$ to detect whether a lane is regular or AV-exclusive lane. If lane $i$ is a regular lane (AV-exclusive lane), \textit{i.e.,} $b_i=0\ (1)$, this sending flow is restricted to the capacity of a regular lane, $\delta q_{lv}$ (AV-exclusive lane, $\delta q_{av}$). {Note that, the candidate AV lanes that are not chosen to be AV-exclusive retains the same characteristics of that of regular lanes with $b_i = 0$.}
\begin{subequations}
	\begin{align}
	\label{send2lv}
	\sum_{k \in A_s}\sum_{j \in \Gamma^+(i)}\Big(\y+\yav\Big) &\leq \delta q_{lv} &&\forall i \in A \setminus A_{av}, \forall t \in T\\
	\label{send3lv}
	\sum_{k \in A_s}\sum_{j \in \Gamma^+(i)}\Big(\y+\yav\Big) &\leq \delta q_{lv}(1-b_i)+\delta q_{av}b_i &&\forall i \in A_{av}, \forall t \in T
	\end{align}
\end{subequations}

{We would like to reiterate that this study assumes a basic level of automation (level 1 automation on the Society of Automotive Engineers automation scale) to be available to the AVs which could be used only on the AV-exclusive lanes. On regular lanes, the AVs are restricted from using these automated facilities and assumed to behave like regular vehicles (LV). Hence, the link capacities of regular lanes are kept fixed as they are assumed to carry vehicles with similar driving behaviour.} 

In the LTM, the receiving flow is defined as the amount of vehicular flow allowed to be received at link $j$ from link $i$ depending on the congestion level and the capacity of link $j$. The receiving flow constraint, as presented in Eq. \eqref{receive1lv}, is derived based on the backward propagation of a congested traffic state from the downstream boundary of a link which reaches the upstream boundary $\frac{L_i}{w_{av}}$ time units later. Here, $w_{av}$ denotes the backward wave speed of the congested traffic state.
{
\begin{multline}
\label{receive1lv}
\sum_{i \in \Gamma^-(j)} \sum_{k \in A_s}  \Big(\y+\yav\Big)\\ \leq \klv L_j-\sum_{k \in A_s}\bigg(\Big(z^{k+}_{j,lv}(t)-z^{k-}_{j,lv}(t_{r,lv})\Big)+\Big(z^{k+}_{j,av}(t)-z^{k-}_{j,av}(t_{r,av})\Big)\bigg)\\\forall j \in A\setminus \{A_r,A_s\}, \forall t \in T \text{ where, }t_{r,lv} = t+\delta-\frac{L_i}{w_{lv}},t_{r,av} = t+\delta-\frac{L_i}{w_{av}}
\end{multline}
}
Similar to Eqs. \eqref{send2lv} and \eqref{send3lv}, Eqs. \eqref{receive2lv} and \eqref{receive3lv} represent the capacity constraint on receiving flow of a link with the binary parameter, $b$, depending on the lane being a regular or candidate AV lane.
\begin{subequations}
	\begin{align}
	\label{receive2lv}
	\sum_{k \in A_s}\sum_{i \in \Gamma^-(j)}\Big(\y+\yav\Big) &\leq \delta q_{lv}\qquad &&\forall j \in A\setminus A_{av},\forall t \in T\\
	\label{receive3lv}
	\sum_{k \in A_s}\sum_{i \in \Gamma^-(j)}\Big(\y+\yav\Big) &\leq \delta q_{lv}(1-b_j)+\delta q_{av}b_j\qquad &&\forall j \in A_{av},\forall t \in T
	\end{align}
\end{subequations}

In the proposed formulation, LVs are restricted from entering an AV-exclusive lane. We formulate integer-linear constraints to implement this restriction in our model as presented in Eqs. \eqref{send2lvav} and \eqref{receive2lvav}. Using the binary lane design variable ($b_i$), the transfer flow of LVs at any timestep is either restricted or kept free, \textit{i.e.}, equal to the capacity of regular lane, for a downstream AV-exclusive or regular lane respectively.
\begin{subequations}
	\begin{align}
	\label{send2lvav}
	\sum_{(o,d) \in W}\sum_{j \in \Gamma^+(i)}\y &\leq (1-b_i)\delta q_{lv} \qquad &&\forall i \in A_{av},\forall t \in T\\
	\label{receive2lvav}
	\sum_{(o,d) \in W}\sum_{i \in \Gamma^-(j)}\y &\leq (1-b_j)\delta q_{lv} \qquad &&\forall j \in A_{av},\forall t \in T
	\end{align}
\end{subequations}

Eqs. \eqref{endtriplv} and \eqref{endtripav} conclude the lane-based LTM formulation ensuring the exit of all the vehicles that entered into the network and reaching their respective destinations at the end of the last timestep ($\bar{t}$).
\begin{subequations}
	\begin{align}
	\label{endtriplv}
	z^{k+}_{k,lv}(\bar{t})&=\sum_{t \in R \setminus {\bar{t}}} D^{o,k}(t)(1-p^{o,k}) &&\forall k \in A_s,\forall (o,k) \in W\\
	\label{endtripav}
	z^{k+}_{k,av}(\bar{t})&=\sum_{t \in R \setminus {\bar{t}}} D^{o,k}(t)p^{o,k} &&\forall k \in A_s,\forall (o,k) \in W
	\end{align}
\end{subequations}

\subsubsection{{OD travel time estimation}}
{The average travel time on each link ($\tau_{lv,l}$, $\tau_{av,l}$) is estimated based on the number of timesteps each vehicle spends on that link, averaged over the total demand of that vehicle class that has travelled through that link as presented in Eqs. \eqref{ttlv1} and \eqref{ttav1}. For each link ($l \in M$), the numerators in Eqs. \eqref{ttlv1} and \eqref{ttav1} represent the total ``system" travel time where ``system" refers to the class of vehicles (LV or AV) travelling on the lanes ($i \in m(l)$) of that link. Here, $m(l)$ denotes the set of lanes in link $l$. The denominators represent the total number of vehicles travelled through that link during the analysis period.
\begin{subequations}
	\begin{align}
	\label{ttlv1}
	\tau_{lv,l} &= \left(\frac{\sum\limits_{i \in m(l)}\sum\limits_{k \in A_s}\sum\limits_{t \in T}\Big(z_{i,lv}^{k+}(t)-z_{i,lv}^{k-}(t)\Big)\delta}{\sum\limits_{i \in m(l)}\sum\limits_{k \in A_s}z_{i,lv}^{k+}(t_n)}\right)&&\forall l \in M \\
	\label{ttav1}
	\tau_{av,l} &= \left(\frac{\sum\limits_{i \in m(l)}\sum\limits_{k \in A_s}\sum\limits_{t \in T}\Big(z_{i,av}^{k+}(t)-z_{i,av}^{k-}(t)\Big)\delta}{\sum\limits_{i \in m(l)}\sum\limits_{k \in A_s}z_{i,av}^{k+}(t_n)} \right)&&\forall l \in M
	\end{align}
\end{subequations}
}

{To estimate the average OD travel time for each vehicle class, we assume that the set of paths, $\Pi^{o,d}$, between each OD pair is known and the travel times of the links ($\tau_{lv,l},\tau_{av,l}$), belonging to these paths are averaged over the number of paths as shown in Eqs. \eqref{lvpath} and \eqref{avpath}.
\begin{subequations}
	\begin{align}
	\label{lvpath}
	\tlv &= \frac{1}{|\Pi^{o,d}|}\sum\limits_{\pi \in \Pi^{o,d}} \sum\limits_{l \in \pi} \tau_{lv,l} && \forall (o,d) \in W\\
	\label{avpath}
	\tav &= \frac{1}{|\Pi^{o,d}|}\sum\limits_{\pi \in \Pi^{o,d}} \sum\limits_{l \in \pi} \tau_{av,l} && \forall (o,d) \in W
	\end{align}
\end{subequations}
}

\subsubsection{{Objective function}}
We choose the total system travel time (TSTT) as the metric of system performance, consisting of travel times of two vehicle classes: LV and AV. The objective function in the proposed formulation minimizes this TSTT as shown in Eq. \eqref{tstt}. 
\begin{align}
\label{tstt}
\min\ (\text{TSTT}_{lv}+\text{TSTT}_{av})
\end{align}

$\text{TSTT}_{lv}$ and $\text{TSTT}_{av}$ can be obtained from the difference in cumulative inflows and outflows of each link, representing the number of LVs/AVs present on that link at each timestep and the number of timesteps they spend on that link. The underlying LTM provides these cumulative inflows and outflows as output as explained in Section \ref{ltm}. 

In the LTM, the cumulative inflows and outflows of each lane at each timestep track the vehicular flow in the network. The difference between these inflows and outflows of a lane at each timestep represents the number of vehicles present in that lane for $\delta$ time units, where $\delta$ is the duration of each timestep. Hence, the sum of these differences over all the lanes, OD pairs and timesteps will provide the TSTT of each vehicle class in the network as follows. 
{
\begin{subequations}
\begin{align}
\label{tsttlv}
\text{TSTT}_{lv} &= \delta \sum_{i \in A \setminus A_s}\sum_{k \in A_s} \sum_{t \in T} \Big(\zp-\zn\Big) \\
\label{tsttav}
\text{TSTT}_{av} &= \delta \sum_{i \in A \setminus A_s}\sum_{k \in A_s} \sum_{t \in T} \Big(\zpav-\znav\Big)
\end{align}
\end{subequations}
}	
We rewrite the objective function of the proposed MINLP presented in Eq. \eqref{tstt} as follows. 
\begin{align}
\label{eqfirst}
\min \delta \sum_{i \in A \setminus A_s}\sum_{(o,d) \in W} \sum_{t \in T} \Big(\zp+\zpav-\zn-\znav\Big)
\end{align}

Note that, as LVs are restricted on AV-exclusive lanes, we fix at least one path with regular lanes between each OD-pair in our model for movement of LVs.

\subsubsection{{Model summary}}
The resulting MINLP formulation $\PF$ represents the freeway network design problem. 
\begin{model}[$\PF$]
\label{mod1}
\begin{braced}
\begin{tabular}{ll}
&$\min$ \text{TSTT} \quad \eqref{eqfirst} \\
&\text{\emph{s.t.:}}\\
&\text{Endogenous demand}\quad \eqref{pav}, \eqref{lvpath}, \eqref{avpath}\\
&\text{Network dynamics}\quad \eqref{demandlv} - \eqref{endtripav}\\
&$\bm{y} \in \Ydom, \bm{z} \in \Zdom, \bm{b} \in \Bdom, \bm{\tau} \in \Tdom, \bm{p} \in \Pdom$
\end{tabular}
\end{braced}
\end{model}

As presented above, $\PF$ involves five sets of variables: transfer flows ($\bm{y}$) with domain $\Ydom = \mathbb{R}_+^{|A_c||\Gamma^+(A_c)||A_s||T|}$, cumulative inflows and outflows ($\bm{z}$) with domain $\Zdom= \mathbb{R}_+^{|A||A_s||T|}$, binary variables for lane allocation ($\bm{b}$) with domain: $\Bdom=\{0,1\}^{|A_{av}|}$, class-wise travel times ($\bm{\tau}$) with domain $\Tdom=\mathbb{R}_+^{2|W|}$, and OD proportion of AVs ($\bm{p}$) domain $\Pdom = {[0,1]^{|W|}}$. Due to the integer variables for lane design ($\bm{b}$) and the nonlinear logit model, $\PF$ may lead to computational tractability issues for bigger networks. In Section \ref{solution}, we deal with these non-linearity issues by introducing a Benders decomposition approach with an embedded fixed-point algorithm and implement it on a freeway network in Section \ref{numerical}.

The outputs of Model $\PF$ can be interpreted as follows. The main output are the lane design variables $\bm{b}$ which indicate which candidate lane should be AV-exclusive in the freeway network. 

The remaining variables are used to account for congestion effects and endogenous travel demand. Travel demand is loaded into the network through the source centroid connectors as expressed in Eqs. \eqref{demlv} and \eqref{demav}. At the completion of the trips, Eqs. \eqref{endtriplv} and \eqref{endtripav} ensures that vehicles leave the network through the sink centroid connectors. On a freeway, these source and sink centroid connectors represent on- and off-ramps respectively. If the network is unable to accept demand due to congestion on the freeway, vehicles may be held at on-ramps which are assumed to have sufficiently large capacities. Since waiting time is penalized in the objective function, the outputs of the proposed formulation can be interpreted as the level of control at the freeway on-ramps.

We analyse the proposed model in the following sections.

\subsection{Fixed-point analysis}
\label{fp}
To motivate the design of a dedicated solution method and to provide insights into the behavior of $\PF$, we consider a simplified version of the model wherein the endogenous demand $\bm{p}$ and the lane design $\bm{b}$ are assumed fixed. This simplified model is called $\SP$ and presented below. 

\begin{model}[$\SP$]
\label{mod2}
\begin{braced}
\begin{tabular}{ll}
&$\min\ \text{TSTT}$ \quad \eqref{eqfirst} \\
&\text{\emph{s.t.:}}\\
&\text{Network dynamics}\quad \eqref{demandlv} - \eqref{endtripav}\\
&$\bm{y} \in \Ydom,\bm{z} \in \Zdom$
\end{tabular}
\end{braced}
\end{model}

The variables involved in Model $\SP$ are transfer flows $\bm{y} \in \Ydom$ and cumulative inflows and outflows $\bm{z} \in \Zdom$. {This model does not involve the average travel times, $\tlv, \tav$, as it represents a parameterized sub-problem that computes a flow pattern under fixed lane design ($b$) and mode split proportion ($\probav$) conditions. Hence, unlike Model \ref{mod1}, in Model \ref{mod2} the lane design and the demand for each vehicle class are fixed. This model facilitates the fixed-point analysis and will be used within the proposed solution method in Section \ref{3.1}}.

\begin{figure}[H]
\includegraphics[clip, trim=22cm 19cm 0.5cm 7cm, width=1.60\textwidth]{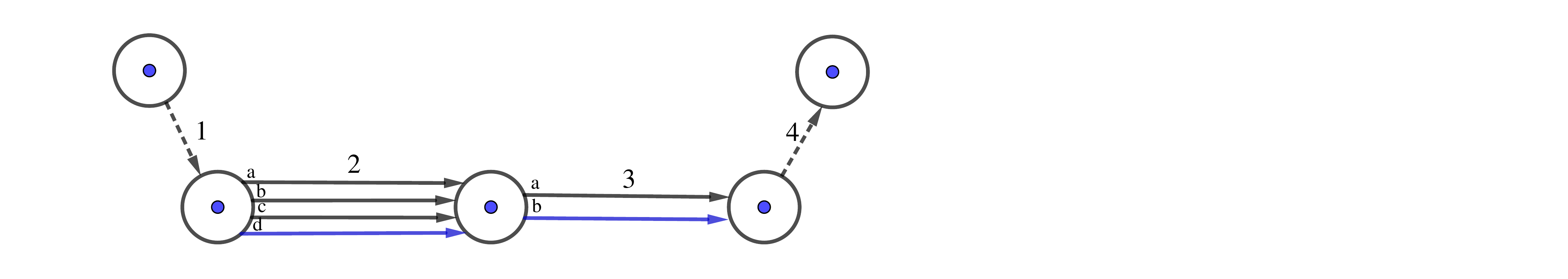}
\caption{Single OD case study network}
\label{fig:singleod}
\end{figure}

We consider a single-OD network with two links with fixed AV-exclusive lanes as illustrated in Figure \ref{fig:singleod}, where the fixed AV lanes are shown in blue. The first link consists of three regular lanes and one AV-exclusive lane followed by a capacity drop on the second link which has one LV and one AV-exclusive lane. We solve $\SP$ for a series of OD proportions of AVs ($\probav$) and calculate the corresponding logit-derived proportion of AVs ($p^{o,d}_{logit}$) using Eq. \eqref{pav} based on the optimal solution of $\SP$. Note that if $p = p_{logit}$, then the demand splits are logit-compatible, i.e. equilibrated, and this solution corresponds to a fixed-point. We define this fixed-point as follows.

\begin{defi}[Fixed-point]
Let $F: \Pdom \rightarrow \Pdom$ be a continuous function of the OD proportion vector $\bm{p} \in \Pdom$. We say that $\bm{p}$ is a fixed-point if $F(\bm{p}) = \bm{p}$. 
\end{defi}

The regular and AV-exclusive lanes differ from each other in terms of capacity and backward wave speed of congestion propagation. In this study network, the capacity (4320 veh/hr) and backward wave speed (28.4 km/hr) of an AV-exclusive lane is taken as double that of an regular lane due to the inter-connectivity of AVs leading to better traffic flow and faster congestion propagation. The length and jam density of the links are 800m and 200 veh/km respectively with a free-flow speed of 90 km/hr for both vehicle classes. The demand is loaded through the source centroid connector 1 into the network. The capacity and jam density of the source and sink centroid connectors are set to very high values with a negligible length for instantaneous loading of demand into the network based on available network capacity. The values of these network parameters are provided in Table \ref{tab:netpar}.

\begin{figure}[!h]
\centering
\includegraphics[width=0.7\textwidth]{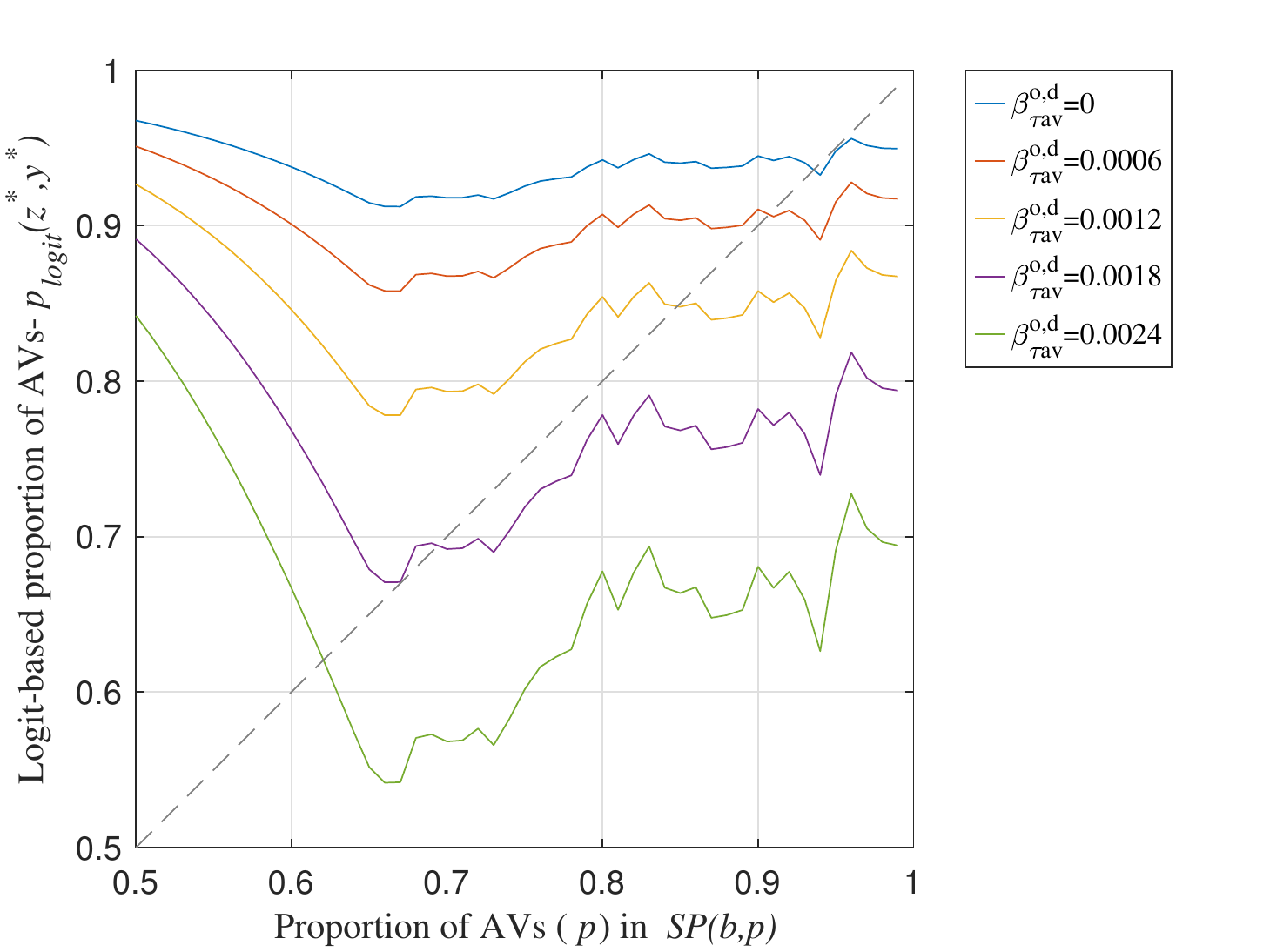}
\caption{{Illustration of the relationship between the proportion of AVs (p) input to SP(b,p) and the logit-based proportion determined based on the solution of SP(b,p) for different $\bav$ values. Each value of $\bav$ corresponds to a different problem setting.}}
\label{fig:ppopt}
\end{figure}

\begin{table}[H]
\caption{Single OD network characteristics}
\label{tab:netpar}
\resizebox{\textwidth}{!}{%
\begin{tabular}{@{}lllllllll@{}}
	\toprule
	Parameters                 & Source (1) & Lane 2a & Lane 2b & Lane 2c & Lane 2d & Lane 3a & Lane 3b & Sink (4) \\ \midrule
	Length (km)                & 0.0001           & 0.8    & 0.8    & 0.8    & 0.8    & 0.8    & 0.8    & 0.0001         \\
	Free-flow speed (km/hr)     & 90               & 90     & 90     & 90     & 90     & 90     & 90     & 90             \\
	Backward wave speed (km/hr) & 12.2               & 12.2     & 12.2     & {12.2}    & 28.4    & 12.2     & 28.4    & 12.2             \\
	Capacity (veh/hr)          & 360000           & 2160   & 2160   & {2160}   & 4320   & 2160   & 4320   & 360000         \\
	Jam density (veh/km)       & 100000           & 200    & 200    & 200    & 200    & 200    & 200    & 100000         \\ \bottomrule
\end{tabular}%
}
\end{table}
We start this experiment by varying the proportion of AVs ($\probav$) from 0.5 to 0.99 in steps of 0.01 and solve Model $\SP$ at each value of $\probav$. For each step, we obtain the average OD travel times of each vehicle class ($\tlv, \tav$) using Eqs. \eqref{lvpath} and \eqref{avpath}. We then calculate $\probav_{logit}$ by substituting $\tlv$ and $\tav$ in Eq. \eqref{pav} along with the coefficients of these travel times ($\blv,\bav$). These coefficients are obtained from a previous study on route choice behaviour of LVs and AVs where the value of time for LV and AV users were found to be \$10 and \$6.5/hr \cite{neeraj2018}. To identify fixed-points ($\probav=\probav_{logit}$), we plot $p$ against $\probav_{logit}$ for different values of $\bav$ while $\blv$ remains fixed, as shown in Figure \ref{fig:ppopt}. 
\begin{figure}[!h]
\centering
\includegraphics[width=0.7\textwidth]{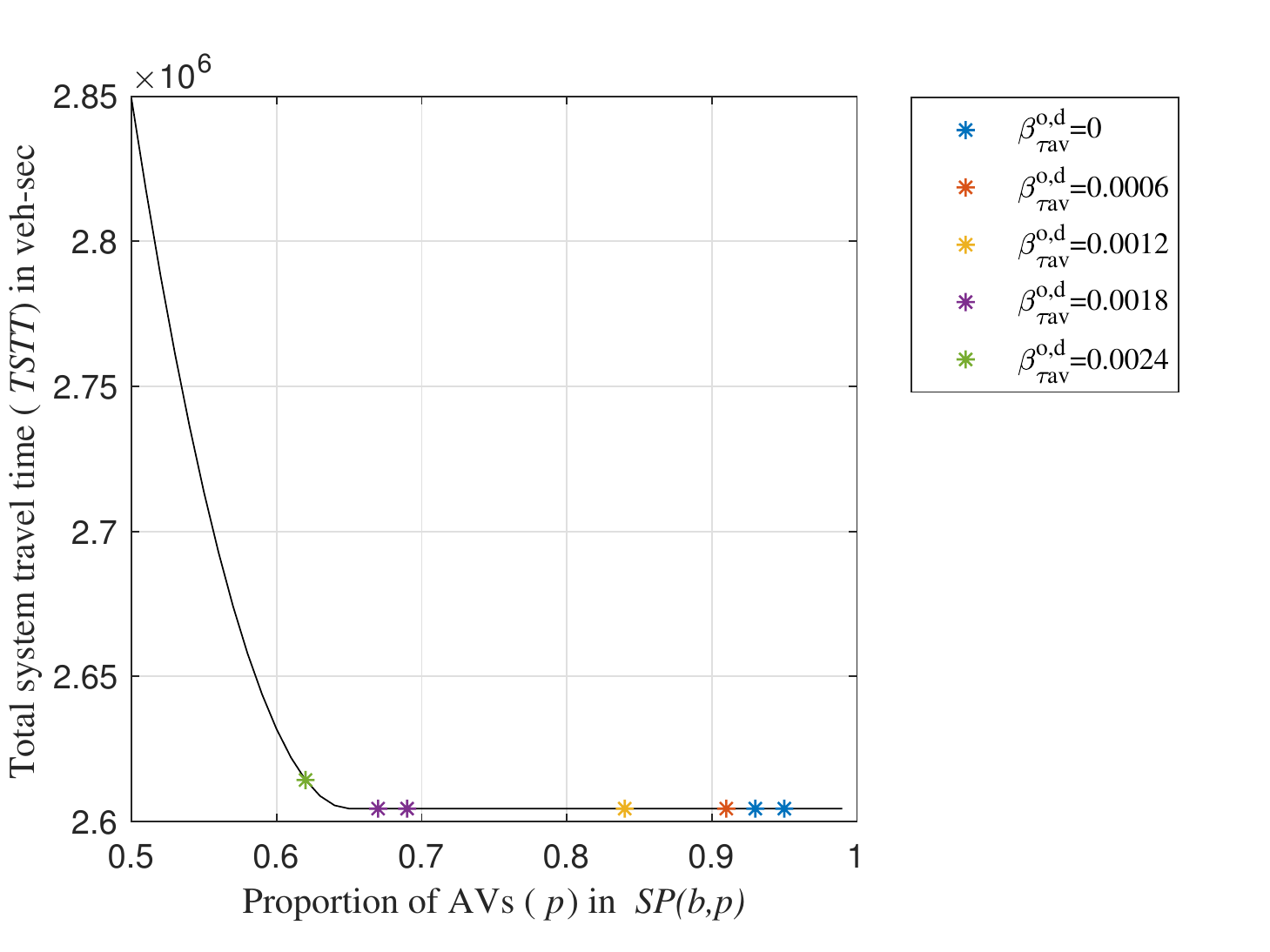}
\caption{Change in TSTT with $p$}
\label{fig:ptstt}
\end{figure}

The dotted line in Figure \ref{fig:ppopt} acts as a reference line to locate fixed-points. Interestingly, the line curve depicting the relationship between $\probav$ and $\probav_{logit}$ in Figure \ref{fig:ppopt} is found to cross this reference line multiple times showing the existence of multiple fixed points in the problem. {Here, each value of $\bav$ corresponds to a different problem setting, and multiple fixed points are observed only for $\bav = 0$ and $\bav = 0.0018$}. Figure \ref{fig:ptstt} depicts the change in the value of the objective function (TSTT) with respect to $p$. {In this Figure, we observe that the points with same colour (blue and purple), representing multiple fixed-points emerging from a same $\bav$, have equal TSTT}. This experiment highlights that fixed-points may be non-unique, but may also correspond to identical TSTT. 

From Figure \ref{fig:ppopt}, we observe that in each case, the line plot depicting the relationship between $p$ and $p_{logit}$ crosses the reference line at least once, referring to the existence of at least one fixed-point satisfying the logit model. We prove this existence of at least one fixed-point with Proposition \ref{prop1}.

\begin{prop}
\label{prop1}
Existence of fixed-point: For a fixed lane design vector $\bm{b} \in \Bdom$, there exists at least one fixed-point such that $F(\bm{p}) = \bm{p}$.

\begin{proof}
We show that we can construct such a continuous function $F(\bm{p})$. Let $h(\bm{p})$ be the function mapping the OD proportion vector $\bm{p}$ to the optimal solution of the linear program $\SP$ for fixed lane design vector $\bm{b}$ as defined by Model \ref{mod2}. Let $\bm{z} \in \Zdom$ be the vector of cumulative inflow and outflow and let $\bm{y} \in \Ydom$ be the vector of transfer flows obtained after solving $\SP$, formally:
\begin{align*}
h&: \Pdom \rightarrow \Zdom \times \Ydom\\
h(\bm{p}) &= (\bm{z}, \bm{y})
\end{align*}
Let $g(\bm{z}, \bm{y})$ be the function mapping the optimal cumulative inflows and outflows, and transfer flows to average, class-based (LV and AV) OD travel times $(\bm{\tau}_{lv},\bm{\tau}_{av}) \in \Tdom$ as defined in Eqs. \eqref{lvpath} and \eqref{avpath}:
\begin{align*}
g&: \Zdom \times \Ydom \rightarrow \Tdom\\
g(\bm{z^\star,y^\star}) &= (\bm{\tau}_{lv},\bm{\tau}_{av})
\end{align*}
Finally, let $f(\bm{\tau}_{lv},\bm{\tau}_{av})$ be the function mapping average class-based OD travel times to OD proportion $\bm{p} \in \Pdom$ via the proposed logit model as defined in \eqref{pav}:
\begin{align*}
f&: \Tdom \rightarrow \Pdom\\
f(\bm{\tau}_{lv},\bm{\tau}_{av}) &= \bm{p}
\end{align*}
Let $F(\bm{p}) = f(g(h(\bm{p})))$, $F$ is a continuous function from the compact convex set $\Pdom$ to itself. By Brouwer's theorem, there exists at least one fixed-point $\bm{p}$ such that $F(\bm{p})=\bm{p}$. 
\end{proof}
\end{prop}  

In Section \ref{solution}, we develop the solution methodology to solve the non-linear Model \ref{mod1} with variable lane design and endogenous demand based on the embedded logit model.

\section{Solution methodology}
\label{solution}
In this section, we propose a decomposition approach with variable lane design and endogenous demand for each vehicle class. The purpose of this development is to identify the crucial locations in a network where providing an AV-exclusive lane will reap the maximum benefit for a given AV demand. As LVs are restricted on AV-exclusive lanes, it is necessary to deploy AV-exclusive lanes judiciously to cater to both vehicle classes keeping the social welfare into perspective. We introduce Benders' decomposition method for this purpose which  iteratively explores possible lane designs in a master problem and, at each iteration, solves a sequence of SODTA problems which is shown to converge to fixed-points representative of logit-compatible demand splits.

\subsection{Benders' decomposition approach}
Benders' decomposition approach eases up the computation burden of a mathematical model by partitioning the overall formulation into a relaxed master problem with mainly integer variables and a subproblem with all the continuous variables. For a detailed review of Benders' approach, one can refer to \citet{RAHMANIANI2017801}. For problems with minimizing objective function, such as the proposed model in this study, the relaxed master problem provides a lower bound at each iteration of Benders' method. Whereas, the subproblem handles the complicated constraints of the original problem which is solved iteratively for each relaxed solution of the master problem.

In our decomposition of Model $\PF$, we only retain binary lane design variables $\bm{b}$ in the master problem $\MP$, which is summarized below:
\begin{model}[$\MP$]
\label{MP}
\begin{braced}
\begin{tabular}{ll}
	&$\min$ $Z$\\
	&\text{\emph{s.t.:}}\\
	&$Z \geq$ \text{Optimality cuts}\\
	&$0 \geq$ \text{Feasibility cuts}\\							
	&$\bm{b} \in \Bdom$, $Z \geq 0$
\end{tabular}
\end{braced}
\end{model}

{In Model \ref{MP}, $Z$ is a surrogate objective function which is an underestimation of the objective function of Model \ref{mod1}, i.e., TSTT. At each iteration until convergence, the Benders' method generates either a feasibility cut or an optimality cut towards obtaining the optimal solution. If the sub-problem solution is feasible, an upper bound is obtained and an optimality cut is derived towards closing the optimality gap between $Z$ and the upper bound. On the other hand, a feasibility cut is generated to eliminate an infeasible solution provided by the subproblem, preventing the master to produce it again.}

The algorithm with Benders' decomposition approach is presented in Algorithm \ref{algo:lbfp}. Here, the subproblem ($\SP$) is initiated by fixing the set of binary lane design variables. Depending on the feasibility of the subproblem, dual prices or rays for each constraint are calculated, followed by solving the master problem which provides the values for the next set of binary lane design variables. This iterative process continues until we reach an exact solution of the objective function.

The subproblem in the proposed formulation consists of two components: the lane-based LTM with system-level objective function and an endogenous demand model. As the endogenous demand model introduces non-linearity in the formulation, we develop a fixed-point algorithm for solving the subproblem as explained in the following subsection. 

\subsection{Fixed-point algorithm}
\label{3.1}
The endogenous demand model is crucial to study the effect of infrastructural changes such as AV-exclusive lanes on AV demand. Note that, the formulation without the endogenous demand model is a useful model by itself as it can estimate the progressive deployment of AV-exclusive lanes in a network corresponding to incremental penetration of AV demand. However, this model does not account for the effect of these AV-exclusive lanes on the demand of each vehicle class.

\begin{algorithm}[!tbp]
	\setstretch{1.0}
	$\bm{p}^0 \gets \frac{1}{e^{\tau_{ff}^{o,d}\big(\beta^{o,d}_{\tau_{av}}-\beta^{o,d}_{\tau_{lv}}\big)}+1}$ (Solve logit model with free-flow travel times)\\
	$\bm{b}^0 \gets \bm{0}$ (No AV-exclusive lanes)\\
	$m \gets 0$\\
	$n \gets 0$\\
	${Z^m \gets 0}$\\
	$UB \gets \infty$\\
	\Repeat{$GAP \leq \epsilon$}{
		$m \gets m + 1$\\
		\Repeat{$\sum_{(o,d) \in W}|p^{o,d,n+1}-p^{o,d,n}|\leq \epsilon_{MSA}$}{
			$n \gets n + 1$\\
			$\bm{y}^n, \bm{z}^n \gets \text{Solve } \texttt{SP}(\bm{b}^m,\bm{p}^n)$\\
			{\For{$l \in M$}{
					$\tau^n_{lv,l} \gets \left(\frac{\sum\limits_{i \in m(l)}\sum\limits_{k \in A_s}\sum\limits_{t \in T}\Big(z_{i,lv}^{k+}(t)-z_{i,lv}^{k-}(t)\Big)\delta}{\sum\limits_{i \in m(l)}\sum\limits_{k \in A_s}z_{i,lv}^{k+}(t_n)}\right)$\\
					$\tau^n_{av,l} \gets \left(\frac{\sum\limits_{i \in m(l)}\sum\limits_{k \in A_s}\sum\limits_{t \in T}\Big(z_{i,av}^{k+}(t)-z_{i,av}^{k-}(t)\Big)\delta}{\sum\limits_{i \in m(l)}\sum\limits_{k \in A_s}z_{i,av}^{k+}(t_n)}\right)$
			}}
			\For{$(o,d) \in W$}{
				{$\tau^{o,d,n}_{lv} \gets \frac{1}{|\Pi^{o,d}|}\sum\limits_{\pi \in \Pi^{o,d}} \sum\limits_{l \in \pi} \tau_{lv,l}$\\
					$\tau^{o,d,n}_{av} \gets \frac{1}{|\Pi^{o,d}|}\sum\limits_{\pi \in \Pi^{o,d}} \sum\limits_{l \in \pi} \tau_{av,l}$} \\
				${\probav} \gets \frac{1}{e^{\big(\beta^{o,d}_{\tau_{av}}\tau^{o,d,n}_{av}-\beta^{o,d}_{\tau_{lv}}\tau^{o,d,n}_{lv}\big)}+1}$\\
				$p^{o,d,n+1} \gets \frac{n}{n+1}p^{o,d,n} + \frac{1}{n+1}{\probav}$
			}
		}
		\If{$\texttt{SP}(\bm{b}^m,\bm{p}^n)$ is infeasible}{
			Generate feasibility cut\\
		}
		\uElse{
			\If{$\text{TSTT}^n<UB$}{
				$UB \gets \text{TSTT}^n$\\
				$\bm{b}^\star \gets \bm{b}^m$\\
				$\bm{\tau}^\star \gets \bm{\tau}^n$\\
				$\bm{p}^\star \gets \bm{p}^n$
			}
			$GAP \gets \frac{UB-{Z^{m-1}}}{UB}$\\
			\If{$GAP>\epsilon$}{
				Generate optimality cut\\
		}}
		$\bm{b}^m,Z^m \gets$ Solve $\MP$
	}
	\Return $UB,\bm{b^\star, p^\star, \tau^\star}$\\
	\caption{Benders decomposition with fixed-point algorithm for \PF}
	\label{algo:lbfp}
\end{algorithm}

We adopt the logit model for endogenously estimate the demand for each vehicle class. As the logit model is nonlinear, we use a fixed-point algorithm to separate the endogenous demand component from the SODTA formulation. The fixed-point algorithm estimates the proportion of AVs ($\probav$) in an iterative process involving the logit model, presented in Eq. \eqref{pav}. The Algorithm \ref{algo:lbfp} proposed in this study keeps the nonlinear logit model out of the MILP, keeping mathematical formulation linear. The new proportion of AV demand is obtained based on the difference in utilities between the modes, substituted in the logit model (as shown in Eq. \eqref{pav}), which is fed back to Model \ref{mod2} for subsequent iterations. We adopt the method of successive averages (MSA) for convergence of the fixed-point algorithm which is based on a predetermined move size along the descent direction. 

\textcolor{blue}{MSA is one of the most widely used solution heuristics in DTA models due to its simplicity \cite{sheffi1985urban,wu2003network,sbayti2007efficient,bell2008attacker}. As MSA depends on predetermined step sizes (decreasing with the iteration index) without the requirement of derivative information for flow-cost mapping, this method does not necessitate finding the optimal move size at every iteration. However, additional conditions regarding continuous first and second derivative of the objective function must be satisfied for MSA algorithms to converge \cite{powell1982convergence}. In this study, we do not provide a formal proof of convergence for MSA. However, our numerical experiments show that the proposed algorithm is able to find feasible solutions for the network design problem at hand.} The iterative process of MSA may disregard the instability in the fixed-point solution. We monitor this instability by checking the value of $\probav$ before and after implementing MSA. 

{The inner repeat loop of Algorithm \ref{algo:lbfp} iterates between two steps: i) solving the subproblem $SP(\bm{b}^m,\bm{p}^n)$ for fixed lane design $\bm{b}^m$ and mode choice proportions $\bm{p}^n$, and ii) adjusting the mode split proportions of each origin-destination (OD) pair $p^{o,d}$ via the mode choice (logit) model.}

{At each iteration of step i), the subproblem is solved with a mode split proportion vector $p^{o,d,n}$ and yields a corresponding flow pattern represented by variables $\bm{z}^n$. This flow pattern is then used in step ii) to used to recalculate the mode split proportions of each OD pair, $p^{o,d}$, via the mode choice (logit) model. The mode split proportions $p^{o,d,n+1}$ are then determined using MSA as convex combination of $p^{o,d}$ and $p^{o,d,n}$.}

{At each iteration of this inner loop, the total absolute value over all OD pairs between the mode choice proportions across two consecutive iterations are compared. The inner repeat loop stops when this total absolute value is less than a predefined threshold ($\epsilon_{MSA}$), hence the inner loop repeat steps i) and ii) until a fixed point is found.}

In the next section, we implement this algorithm on a single-OD and a multi-OD network.

\section{Numerical experiments}
\label{numerical}
The proposed formulation is implemented on the same single OD network (Figure \ref{fig:singleod}) presented in Section \ref{problem} along with a multi-OD freeway network (Figure \ref{fig:bignet}). 

\subsection{Single OD freeway network with fixed lane design}
\begin{figure}[H]
\centering
\centering
\includegraphics[width=0.8\linewidth]{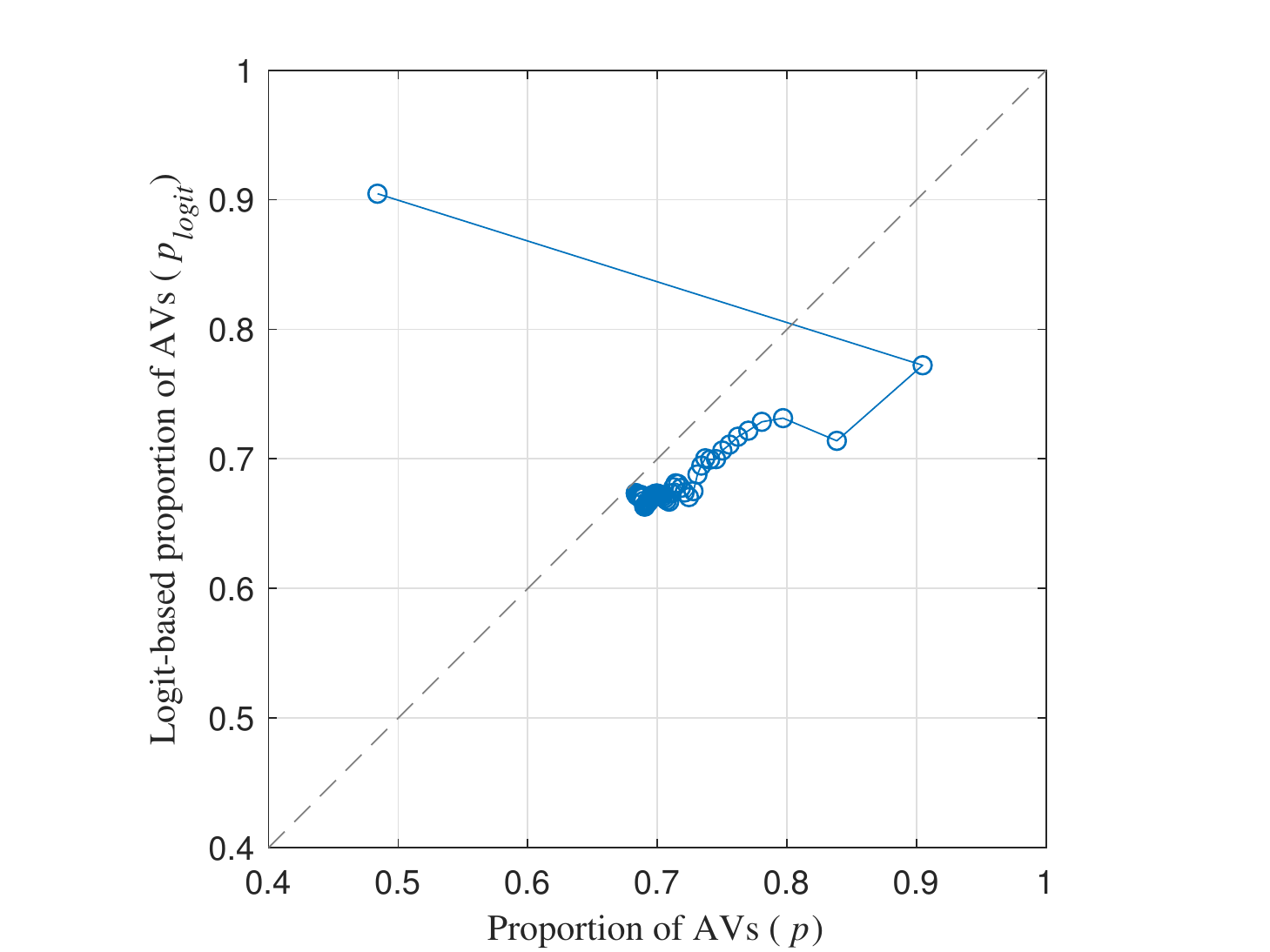}
\captionsetup{justification=centering}
\captionof{figure}{Convergence of FP algorithm \\(with fixed lane design)}
\label{fig:m2ppl}
\end{figure}

We implement the solution methodology presented in Section \ref{solution} on the single OD freeway network with fixed lane design for AVs. As mentioned earlier, the presence of the binary lane design variable makes the proposed formulation a mixed-integer linear program, whereas, the endogenous demand model introduces non-linearity due to the structure of the logit model. We develop the fixed-point algorithm to circumvent non-linearity and Benders decomposition method to handle the binary lane design variable in the master problem. To understand the performance of the fixed-point algorithm alone, we first implement Model \ref{mod2} on the single OD network with fixed dedicated lanes for AVs. We initialize the algorithm with a proportion of AV obtained from Eq. \eqref{pav} substituting the free-flow travel time between the OD pairs. Model \ref{mod2} is solved at each iteration of this algorithm. The value of $\probav$ is updated at each iteration with the proportion ($p_{logit}$) obtained by the endogenous demand model until convergence where a fixed-point is reached ($p = p_{logit}$). The network characteristics of this test network is same as presented in Table \ref{tab:netpar}. A time-varying demand profile is selected for this analysis which is loaded into the network every 2 minutes for the first 8 minutes of a total analysis period of 100 minutes. The total demand is 2950 vehicles which includes both LVs and AVs. Note that, the objective function of the proposed formulation takes into account the waiting time of vehicles at the source centroid connector depending on the available capacity in the network. {We consider a timestep ($\delta$) of 30 seconds, satisfying the Courant-Friedrichs-Lewy (CFL) condition \citep{courant1928partiellen} of LTM which states that this $\delta$ is required to be at most equal to the free-flow travel time ($\delta \leq \frac{L_{min}}{v_f}$) of the shortest link in the network to prevent vehicles from moving out of a link during the time interval.} At each timestep, the cumulative inflows and outflows of each link are updated by the underlying LTM. For a single OD network with fixed dedicated lanes for AVs with $\bav = 0.0018$, we plot this convergence of the fixed-point algorithm in Figure \ref{fig:m2ppl}.

Figure \ref{fig:m2ppl} shows that the algorithm converges to a fixed-point after 67 iterations at $p = 0.68$ which is one of the fixed points presented by the purple dot in Figure \ref{fig:ptstt}. At convergence, the value of the objective function (TSTT) is found to be 2604430 veh-sec with a computation time of 24.067 seconds.

In this case study on a single OD network with fixed AV-exclusive lanes, we observe that the fixed-point algorithm performs well with fast convergence. In Subsection \ref{sub:4.2}, we introduce the Benders' decomposition method to handle variable lane design problem on a multi-OD network. 

\subsection{Multi-OD single path freeway network}
\label{sub:4.2}
\begin{figure}[H]
\includegraphics[trim={0 30cm 0 30cm},clip,width=\linewidth]{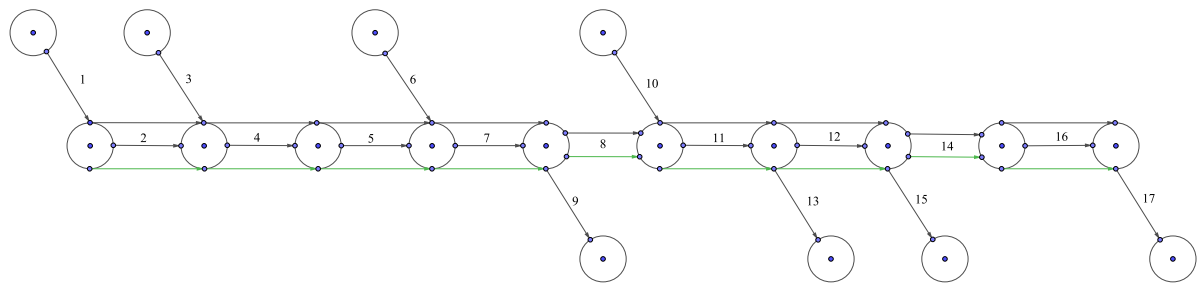}
\caption{Multi-OD network}
\label{fig:bignet}
\end{figure}

We implement the proposed model with Benders method and fixed-point algorithm on a multi-OD freeway network, presented in Figure \ref{fig:bignet}. This 27km long freeway consists of 6 OD pairs, 9 links, 4 on-ramps and 4 off-ramps. Each of these links has either 2 or 3 lanes as shown in Figure \ref{fig:bignet} where one lane in each link, showed in green colour, belongs to the set of candidate AV lanes ($A_{av}$) which could potentially converted into an AV-exclusive lane for better system performance. The total vehicular demand is considered to be 3000 with an analysis period of 50 mins where the demand is loaded into the network through the 4 on-ramps over a period of first 10 mins of the analysis. {This total demand is selected in such a way that it creates significant congestion in the network but still ensures that all vehicles reach their destinations within the analysis period.}The free-flow speed (90 km/hr), capacities (2160 and 4320 veh/hr for regular and AV-exclusive lanes respectively) and the backward wave speeds (12.2 and 28.4 km/hr for regular and AV-exclusive lanes respectively) are considered the same as the single OD network (Figure \ref{fig:singleod}) to match the fundamental diagram of traffic flow, showed in Figure \ref{fig:fd}. The traffic flow propagation is captured every minute based on the lane-based LTM.

We design 7 experiments to study the performance of our algorithm. Apart from the base case, the proposed algorithm is implemented for different demand ($\pm25\%$), capacity of AV lanes ($\pm25\%$) and coefficient of average AV travel times ($\pm25\%$). The performance of Algorithm \ref{algo:lbfp} for these different scenarios is presented in Table \ref{tab:performance2}. 

\begin{table}[H]
\caption{Performance of Algorithm 1 on a multi-OD freeway network for different scenarios}
\label{tab:performance2}
\resizebox{\textwidth}{!}{%
\begin{tabular}{@{}llllllll@{}}
	\toprule
	Parameter                                                            & Base    & -25\% Demand & +25\% Demand & -25\% $q_{av}$ & +25\% $q_{av}$ & -25\% $\beta_{av}$ & +25\% $\beta_{av}$ \\ \midrule
	\begin{tabular}[c]{@{}l@{}}Nb of converted AV    lanes\end{tabular} & 7(9)    & 7(9)         & 6(9)         & 5(9)           & 8(9)           & 7(9)               & 7(9)               \\
	CPU time (mins)                                                      & 119.88  & 71.57        & 144.52       & 26.51          & 255.94         & 99.98              & 98.45              \\
	Nb of Benders ($m$)                                                        & 55      & 64           & 40           & 14             & 144            & 53                 & 49                 \\
	Nb of FP iterations ($n$)                                                            & 269     & 166          & 395          & 73             & 817            & 261                & 237                \\
	Nb of FP per Benders ($n/m$)                                                            & 4.89     & 2.59          & 9.88          & 5.21             & 5.67            & 4.92                & 4.84                \\
	TSTT (veh-sec)                                                       & 2248830 & 1452190      & 3255250      & 2316410        & 2216510        & 2247660            & 2250430            \\
	\%reduction in TSTT                                                  & 11.43   & 9.41         & 12.16        & 8.77           & 12.71          & 11.48              & 11.37              \\ \bottomrule
\end{tabular}%
}
\end{table}

We observe that it is not beneficial to convert all candidate AV lanes to AV-exclusive lanes. {The output from the sensitivity analysis on demand shows that 7 out of 9 candidate lanes are converted to AV-exclusive lanes for the base case. This remains the same for the reduced demand case as well. However, for an increased demand, a lesser number of lanes are converted to AV lanes (6). This could be due to increased the LV demand which requires more regular lanes in the network.} In case of an increased AV lane capacity, the system performance is significantly enhanced by allocating more AV-exclusive lanes (8) compared to the reduced AV-lane capacity case (5). On the other hand, the coefficient of average AV travel times is found to have no effect on AV lane design in the network. These deployed AV-exclusive lanes are found to decrease the TSTT by around 10\% in all the scenarios, with a maximum improvement of 12.71\% for increased AV-lane capacity. Hence, AV-exclusive lanes are always found to have a positive impact in network performance.

The computation time of the proposed algorithm varied from around 0.5 to 4 hours on an Intel(R) Core(TM) i7-6700 @3.40GHz CPU with 16GB RAM. The base case was found to converge within 2 hours with 55 iterations of Benders and a total of 269 fixed-point iterations. Interestingly, the change in the capacity of AV-exclusive lanes is found to have the maximum effect on computation time of the proposed algorithm with a 10 times increment in computation time for inflated capacity compared to the deflated capacity of AV lanes. Similar increment in the number of Benders and fixed-point iterations is also observed in these cases with inflated and deflated capacities.  

\begin{table}[H]
\caption{Optimal proportion of AVs for each OD pair for different scenarios}
\label{tab:pavs}
\resizebox{\textwidth}{!}{%
\begin{tabular}{@{}lllllllll@{}}
\toprule
\multirow{2}{*}{OD pairs} & \multirow{2}{*}{Initial $\probav$} & \multicolumn{7}{c}{Optimal $\probav$}                                                                          \\ \cmidrule(l){3-9} 
&                                    & Base & -25\% Demand & +25\% Demand & -25\% $q_{av}$ & +25\% $q_{av}$ & -25\% $\beta_{av}$ & +25\% $\beta_{av}$ \\ \midrule
(1,13) & 0.33 & 0.73 & 0.70 & 0.76 & 0.74 & 0.72 & 0.82 & 0.62 \\
(3,9) & 0.41 & 0.63 & 0.63 & 0.62 & 0.64 & 0.63 & 0.68 & 0.57 \\
(3,17) & 0.28 & 0.76 & 0.75 & 0.75 & 0.76 & 0.75 & 0.83 & 0.64 \\
(6,15) & 0.38 & 0.66 & 0.65 & 0.70 & 0.66 & 0.67 & 0.74 & 0.57 \\
(10,15) & 0.44 & 0.65 & 0.58 & 0.68 & 0.63 & 0.64 & 0.69 & 0.61 \\
(10,17) & 0.38 & 0.70 & 0.64 & 0.73 & 0.68 & 0.70 & 0.76 & 0.64              \\ \bottomrule
\end{tabular}%
}
\end{table}

Table \ref{tab:pavs} summarizes the optimal proportion of AVs ($\probav$) obtained from the endogenous demand model in different scenarios. We initialize the values of $\probav$ by substituting free-flow travel times of vehicles in the logit model Eq. \eqref{probav}. Due to the structure of the logit model, the initial $\probav$ values are found to be inversely proportional to the free-flow travel time between OD pairs and vary from 0.28 to 0.44. As the free-flow travel times on regular and AV-exclusive lanes are the same, this variation is solely due to the difference in the coefficient of travel times in the logit model.
{Since the FNDP is a nonconvex problem, and since the proposed solution method cannot guarantee a global optima, it is possible that Algorithm \ref{algo:lbfp} converge to a different solution if the initial conditions change. However, in our numerical experiments, we found that changing the initial conditions did not change the solution returned by Algorithm \ref{algo:lbfp}. However, we observe a faster convergence of the proposed algorithm with the free-flow travel time based initial $\probav$ values compared to a lower starting value (e.g., $\probav$ = 0).}.

In congested traffic conditions the proposed algorithm estimates these $\probav$ to be around 0.57 to 0.83. In the base case, the OD pair (3,17) is found to have the maximum value of $\probav$ with a value of 0.76, whereas the OD pair (3,9) requires the least proportion with a value of 0.63. With a reduction in total demand, this proportion is found to be reduced by a maximum of 10\%. Whereas, with increased demand, the maximum increment in $\probav$ is only around 4.6\%. This could be due to the increased total demand which require fewer AV-exclusive lanes, as shown in Table \ref{tab:performance2}, to cater for the increased LV demand. Hence, increasing $\probav$ may not provide much improvement in system level performance. On the other hand, $\probav$ values are not found to be affected significantly with inflated and deflated capacity of AV-exclusive lanes. However, the different coefficients of the average travel times of AVs ($\beta_{av}$) are found to considerably affect the AV proportions. A reduction in this coefficient means that average travel time will create less disutility for AVs rendering this mode more attractive than AVs and vice-versa. The output from the proposed algorithm shows similar effect with a maximum increment of 12\% in $\probav$ for a 25\% reduction in $\beta_{av}$ and a maximum of 15.7\% decrease with an increment in $\beta_{av}$.

We present the average travel times LVs and AVs on each link for all the 7 experiments in Tables \ref{tab:ttlv} and \ref{tab:ttav}. In these tables, the highlighted rows are source centroid connectors where vehicles queue up to enter into the network. As mentioned earlier, the proposed algorithm triggers 7 out 9 candidate AV lanes into AV-exclusive lanes. These exclusive lanes are found to have a significant effect on the average travel times of LVs with a maximum reduction of 45\% for the base case. Similar improvement is observed for other demand scenarios as well with a maximum reduction of 41 and 46\% respectively for LVs. These maximum reductions are mainly observed on the waiting times of vehicles on on-ramps. With a deflated capacity of the proposed AV lanes, the average reduction of average travel times on congested links is found to be around 8\% with a maximum reduction of 23\% in the first on-ramp of the network. Whereas, the inflated capacity provided a average reduction of 12\% in the travel times with a maximum reduction of 52\%. In the case with reduced $\bav$ value, the average reduction of travel times for LVs are found to be 5.6\% compared to the case with increased $\bav$ value where this reduction is 14\%. This output is quite intuitive as the reduced a $\bav$ will create less impedance on the mode choice decisions of users based on AV travel times, leading to more AVs in the network and less travel time benefits of for LV users.         

\begin{table}[H]
\centering
\caption{Changes in the average link travel times (in seconds) of LVs in different traffic scenarios with  highlighted source centroid connectors}
\label{tab:ttlv}
\resizebox{\textwidth}{!}{%
\begin{tabular}{lllllllllll}
	\toprule
\multirow{2}{*}{Links} & \multicolumn{2}{c}{Base} & \multicolumn{2}{c}{-25\% Demand} & \multicolumn{2}{c}{+25\% Demand} & -25\% $q_{av}$ & +25\% $q_{av}$ & -25\% $\beta_{av}$ & +25\% $\beta_{av}$ \\ \cmidrule(l){2-3} \cmidrule(l){4-5} \cmidrule(l){6-7} \cmidrule(l){8-8} \cmidrule(l){9-9} \cmidrule(l){10-10} \cmidrule(l){11-11} 
	& Initial & Optimal & Initial & Optimal & Initial & Optimal & Optimal & Optimal & Optimal & Optimal \\ \midrule
\rowcolor{yellow}
1 & 277.88 & 174.828 & 151.166 & 106.739 & 406.912 & 240.976 & 212.511 & 131.875 & 131.725 & 183.033 \\
2 & 181.729 & 185.284 & 146.044 & 120.791 & 220.085 & 262.399 & 193.588 & 178.963 & 240.452 & 170.738 \\
\rowcolor{yellow}
3 & 306.436 & 165.757 & 205.679 & 120.054 & 409.129 & 219.676 & 193.067 & 166.291 & 147.396 & 197.611 \\
4 & 140.619 & 137.943 & 130.73 & 121.139 & 147.593 & 152.443 & 160.908 & 139.039 & 149.653 & 139.45 \\
5 & 147.981 & 142.416 & 124.657 & 127.542 & 141.599 & 130.418 & 139.193 & 140.259 & 131.899 & 133.287 \\
\rowcolor{yellow}
6 & 169.538 & 147.117 & 120.405 & 109.776 & 317.411 & 294.197 & 133.437 & 191.978 & 207.007 & 119.886 \\
7 & 134.415 & 129.13 & 121.8 & 122.169 & 153.188 & 124.741 & 123.178 & 124.885 & 127.869 & 122.729 \\
8 & 131.674 & 121.193 & 126.358 & 121.159 & 128.367 & 128.853 & 125.93 & 121.148 & 120 & 120 \\
\rowcolor{yellow}
10 & 210.512 & 235.873 & 144.859 & 84.8162 & 384.315 & 352.003 & 225.173 & 205.47 & 247.179 & 227.634 \\
11 & 120 & 120 & 120 & 120 & 120 & 120 & 120 & 120 & 120 & 120 \\
12 & 120 & 120 & 120 & 120 & 120 & 120 & 120 & 120 & 120 & 120 \\
14 & 120 & 120 & 120 & 120 & 120 & 120 & 120 & 120 & 120 & 120 \\
16 & 120 & 120 & 120 & 120 & 120 & 120 & 120 & 120 & 120 & 120 \\ \hline
\end{tabular}%
}
\end{table}

With the proposed AV-exclusive lanes, the average travel times of AVs are found to be reduced by an average of 15\%. Note that, these average travel times correspond to AVs on both regular and exclusive lanes. The maximum reduction of 60\% is observed for the waiting time on source connector 10. The average reductions in travel times are around 16.4\% and 13.4\% for the reduced and increased demand cases respectively. The deflated and inflated capacities of the AV-exclusive lanes reduced the average AV travel times of AVs by 12.8 and 17.5\% respectively with a maximum reduction of 64\% for the inflated capacity case. With an decreased and increased $\bav$ the average reduction in the average travel times are around 16.4 and 15\% respectively.

\begin{table}[H]
\centering
\caption{Changes in the average link travel times (in seconds) of AVs in different traffic scenarios with highlighted source centroid connectors}
\label{tab:ttav}
\resizebox{\textwidth}{!}{%
\begin{tabular}{@{}lllllllllll@{}}
	\toprule
	\multirow{2}{*}{Links} & \multicolumn{2}{c}{Base} & \multicolumn{2}{c}{-25\% Demand} & \multicolumn{2}{c}{+25\% Demand} & -25\% $q_{av}$ & +25\% $q_{av}$ & -25\% $\beta_{av}$ & +25\% $\beta_{av}$ \\ \cmidrule(l){2-3} \cmidrule(l){4-5} \cmidrule(l){6-7} \cmidrule(l){8-8} \cmidrule(l){9-9} \cmidrule(l){10-10} \cmidrule(l){11-11}   
	& Initial & Optimal & Initial & Optimal & Initial & Optimal & Optimal & Optimal & Optimal & Optimal \\ \midrule
\rowcolor{yellow}
1 & 228.002 & 194.321 & 107.748 & 83.7857 & 386.231 & 279.905 & 216.835 & 183.595 & 173.859 & 195.476 \\
2 & 174.544 & 174.023 & 146.698 & 139.789 & 211.304 & 201.096 & 190.852 & 162.101 & 174.372 & 169.77 \\
\rowcolor{yellow}
3 & 244.744 & 202.592 & 165.435 & 110.367 & 362.77 & 309.514 & 213.117 & 188.538 & 192.961 & 205.334 \\
4 & 144.699 & 139.908 & 128.386 & 124.517 & 166.232 & 144.193 & 152.298 & 132.189 & 148.148 & 140.269 \\
5 & 132.483 & 128.633 & 131.94 & 126.396 & 140.082 & 136.263 & 120.497 & 140.006 & 130.037 & 138.55 \\
\rowcolor{yellow}
6 & 189.882 & 147.155 & 145.978 & 95.2195 & 240.574 & 254.529 & 105.3 & 145.076 & 145.994 & 146.61 \\
7 & 148.247 & 126.342 & 127.908 & 120 & 167.794 & 126.85 & 123.784 & 133.235 & 124.268 & 126.927 \\
8 & 134.795 & 127.385 & 121.697 & 122.503 & 159.399 & 139.126 & 136.804 & 124.822 & 122.682 & 123.774 \\
\rowcolor{yellow}
10 & 331.996 & 138.413 & 128.935 & 77.064 & 353.871 & 260.086 & 185.354 & 118.813 & 149.424 & 138.433 \\
11 & 120 & 120 & 120 & 120.245 & 120 & 120.416 & 120.264 & 120.408 & 120.219 & 120 \\
12 & 120 & 120 & 120 & 120 & 120 & 120 & 120 & 120 & 120 & 120 \\
14 & 120 & 120 & 120 & 120 & 120 & 120 & 120 & 120 & 120 & 120 \\
16 & 120 & 120 & 120 & 120 & 120 & 120 & 120 & 120 & 120 & 120 \\ \bottomrule
\end{tabular}%
}
\end{table}

{To analyse the performance of the model with more number of candidate AV lanes in a network, another experiment is designed with two candidate lanes per link. 
With single candidate lane on each link, the previous experiment has 9 binary lane allocation variables. This number would increase to 16 for the same network with two candidate lanes per link. Based on our analysis on the proposed algorithm, the number of Benders iterations is found to be quite sensitive with the number of binary lane allocation variables with more variables leading to longer computation time as shown in Table \ref{16bin}. In this comparative analysis, the computation time for the experiment with 16 candidate AV lanes is restricted to 6 hours. At the end of 6 hours, the algorithm activates 13 out of 16 AV-exclusive lanes after 293 and 1143 Benders and fixed-point iterations respectively with 14.94\% reduction in TSTT compared to a network with no AV-exclusive lanes.}

\begin{table}[H]
	\centering
	\caption{Comparative analysis on number of candidate AV lanes}
	\label{16bin}
	\begin{tabular}{lll}
		\hline
		Parameter                    & 9 candidate AV lanes     & 16 candidate AV lanes \\ \hline
		Nb of converted AV   lanes   & 7(9)						& 13(16)                                      \\
		CPU time (hours)             & 2.00                    & 6.02                                 \\
		Nb of Benders ($m$)          & 55                       & 293                                         \\
		Nb of FP iterations ($n$)    & 269                      & 1143                                        \\
		Nb of FP per Benders ($n/m$) & 4.89                     & 3.90                                 \\
		TSTT (veh-sec)               & 2248830                  & 2159760                                     \\
		\%reduction in TSTT          & 11.43              & 14.94                                 \\ \hline
	\end{tabular}%
\end{table}

\subsection{{Multi-OD multi-path freeway network}}

\begin{figure}[H]
	\includegraphics[trim={0 30cm 0 30cm},clip,width=\linewidth]{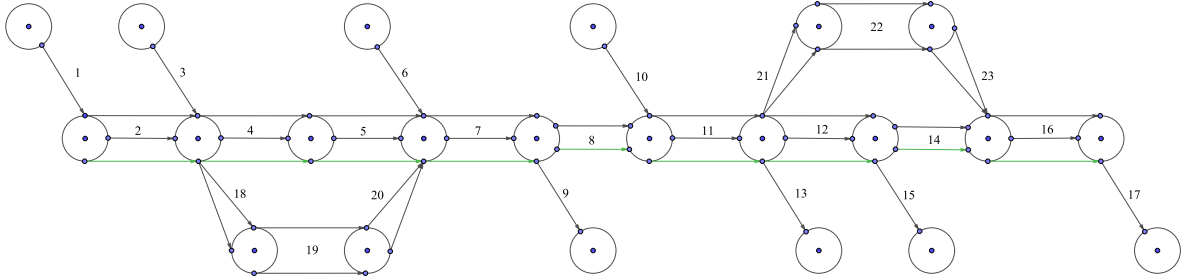}
	\caption{Multi-OD multi-path network}
	\label{fig:multipath}
\end{figure}

{
	In this section, we implement the proposed algorithm on another freeway network with a similar topology but six additional links which provide multiple paths between each OD pair. As the algorithm encapsulates a system optimum dynamic traffic assignment model, the route choice decisions of the users are governed by the system with an objective to minimize overall travel time in the network. The additional links, 18, 19, 20, 21, 22 and 23, are shown in Figure \ref{fig:multipath} where each link consists of two lanes with no candidate AV lanes. In the base case, the length of each link is considered to be 2 km with each lane capacity being 2160 veh/hr and a free-flow speed of 90 km/hr. Along with the base case, four additional experiments are carried out to understand the sensitivity of the proposed algorithm on the length and capacity of the newly added detour links.  
}

{
	In the base case, the algorithm converged in 57 minutes after 11 Benders and 57 fixed-point iterations respectively, converting 4 out of 9 candidate AV lanes into AV-exclusive lanes as presented in Table \ref{tab:performance}. Although, the demand remained the same as the previous experiment with single path between OD pairs, the number of converted AV lanes are found to be lesser in this network with multiple paths. This might be due to the existence of multiple paths which reduces the necessity of high number of AV-exclusive lanes for an overall improvement in the system performance. In the next two iterations, the length of these additional links are varied creating longer and shorter detour options for the users. In the case with longer detour and shorter detours, the length of each of these links are increased to 3 km and decreased to 1.5 km respectively. It is observed that the model converts 7 out of 9 candidate AV lanes into AV-exclusive lanes for the longer detour scenario compared to 4 AV-exclusive lanes in the case with shorter detour. This might be due to longer detour paths not being attractive in terms of travel time triggering more AV lanes on the main road and vice-versa in case of shorter detour paths.}

{
	Two additional experiments involve variation in the capacities of these detour links. It is observed that with decreased capacity, more number of exclusive lanes (4) are activated compared to the case with increased capacity (3). With decreased capacity of the detour links, more vehicles might have chosen the main road and vice-versa leading to such lane allocation outcome. However, almost all the 5 scenarios are observed to improve the TSTT by around 11\% from a setup with no exclusive lanes. 
}

\begin{table}[H]
	\caption{Performance of Algorithm 1 on the multi-OD multi-path freeway network}
	\label{tab:performance}
	\resizebox{\textwidth}{!}{%
		\begin{tabular}{llllll}
			\hline
			Parameter                    & Base    & Longer   detour & Shorter detour & Decreased   Capacity & Increased   Capacity \\ \hline
			Nb of converted AV   lanes   & 4(9)    & 7(9)            & 4(9)           & 4(9)                 & 3(9)                 \\
			CPU time (mins)              & 57.42   & 119.42  & 70.86   & 50.44   & 58.03   \\
			Nb of Benders ($m$)          & 11      & 28      & 13      & 10      & 11      \\
			Nb of FP iterations ($n$)    & 57      & 402     & 143     & 40      & 81      \\
			Nb of FP per Benders ($n/m$) & 5.18    & 14.36   & 4.25    & 4.00    & 7.36    \\
			TSTT (veh-sec)               & 2244070 & 2248600 & 2163500 & 2246600 & 2242820 \\
			\%reduction in TSTT          & 11.13   & 11.44   & 11.74   & 11.28   & 10.93   \\ \hline
		\end{tabular}%
	}
\end{table}

Table \ref{tab:pavsnew} presents the optimal proportion of AVs ($\probav$) for each OD pair for different scenarios. For the base case, the optimal $\probav$ is found to vary between 0.6 to 0.74, whereas, for longer and shorter detour experiments, it varied from 0.59 to 0.77. The change in the optimal $\probav$ is not found to be significant across different experiments for an OD-pair. 

\begin{table}[H]
	\caption{Optimal proportion of AVs for each OD pair for different scenarios}
	\label{tab:pavsnew}
	\resizebox{\textwidth}{!}{%
		\begin{tabular}{lllllll}
			\hline
			\multirow{2}{*}{OD pairs} & \multirow{2}{*}{Initial $\probav$} & \multicolumn{5}{c}{Optimal $\probav$}                                               \\ \cline{3-7} 
			&                                    & Base & Longer detour & Shorter detour & Decreased   Capacity & Increased   Capacity \\ \hline
			(1,13)  & 0.33 & 0.74 & 0.75 & 0.73 & 0.71 & 0.74 \\
			(3,9)   & 0.41 & 0.62 & 0.64 & 0.63 & 0.62 & 0.64 \\
			(3,17)  & 0.28 & 0.74 & 0.77 & 0.74 & 0.75 & 0.76 \\
			(6,15)  & 0.38 & 0.66 & 0.68 & 0.66 & 0.65 & 0.66 \\
			(10,15) & 0.44 & 0.63 & 0.68 & 0.64 & 0.64 & 0.65 \\
			(10,17) & 0.38 & 0.60 & 0.64 & 0.59 & 0.60 & 0.61 \\ \hline
		\end{tabular}%
	}
\end{table}

Tables \ref{tab:lvttmultipath} and \ref{tab:avttmultipath} present the average link travel times of LVs and AVs respectively for different scenarios. In these tables, the source centroid connectors and the detour links are highlighted in yellow and cyan respectively. In the base case, the provision of exclusive lanes is found to reduce the link travel times by a maximum of 48\%. For longer and shorter detour experiments, the maximum reductions are found to be around 23\% whereas for the decreased and increased capacity cases these values are 50 and 25\% respectively. However, for the source centroid connector 6, the AV-exclusive lanes are found to increase the waiting time by 56 and 43\% respectively while minimising overall travel time in the network.

The reduction in the average link travel times is found to be more prominent for AVs as presented in Table \ref{tab:avttmultipath}. For the base case, the maximum reduction in travel time occurred for the source centroid connector 6 with 62\%. For longer and shorter detour lengths, this reduction is found to be 66 and 77\% respectively. This massive reduction in average travel times shows the extent of benefits of AV-exclusive lanes in a network.	    
\begin{table}[!tbp]
	\caption{Changes in the average link travel times (in seconds) of LVs in different traffic scenarios; (source centroid connectors and detour links are highlighted in yellow and cyan respectively)}
	\label{tab:lvttmultipath}
	\resizebox{\textwidth}{!}{%
		\begin{tabular}{lllllllllll}
			\hline
			\multirow{2}{*}{Links} & \multicolumn{2}{c}{Base} & \multicolumn{2}{c}{Longer detour} & \multicolumn{2}{c}{Shorter detour} & \multicolumn{2}{c}{Decreased capacity} & \multicolumn{2}{c}{Increased capacity} \\ \cmidrule(l){2-3} \cmidrule(l){4-5} \cmidrule(l){6-7} \cmidrule(l){8-8} \cmidrule(l){9-9} \cmidrule(l){10-10} \cmidrule(l){11-11} 
			& Initial     & Optimal    & Initial         & Optimal        & Initial         & Optimal        & Initial        & Optimal        & Initial            & Optimal            \\ \midrule
			\rowcolor{yellow}
			1     & 224.84  & 188.21  & 113.55  & 117.09  & 111.07  & 199.47  & 227.14  & 148.63  & 162.16  & 120.31  \\
			2                                       & 147.44  & 136.01  & 292.83  & 222.72  & 158.30  & 160.23  & 221.77  & 181.45  & 180.46  & 146.61  \\
			\rowcolor{yellow}
			3     & 196.71  & 102.22  & 83.35   & 95.86   & 166.69  & 148.32  & 265.07  & 130.74  & 153.80  & 114.38  \\
			4                                       & 183.51  & 160.16  & 195.82  & 157.88  & 145.58  & 159.84  & 167.17  & 150.32  & 206.65  & 163.00  \\
			5                                       & 148.10  & 146.51  & 137.73  & 139.30  & 142.90  & 166.74  & 155.12  & 144.55  & 177.45  & 168.55  \\
			\rowcolor{yellow}
			6     & 234.02  & 168.68  & 119.05  & 185.99  & 166.54  & 180.08  & 106.24  & 152.26  & 207.40  & 159.39  \\
			7                                       & 132.99  & 121.48  & 143.13  & 123.16  & 126.43  & 121.17  & 140.29  & 122.94  & 139.26  & 121.58  \\
			8                                       & 128.15  & 120.00  & 125.15  & 120.00  & 127.39  & 120.00  & 135.20  & 120.00  & 129.48  & 120.00  \\
			\rowcolor{yellow}
			10    & 214.76  & 207.15  & 322.91  & 268.06  & 172.38  & 215.12  & 216.13  & 216.03  & 260.33  & 195.38  \\
			11                                      & 120.00  & 120.00  & 120.00  & 120.00  & 120.00  & 120.00  & 120.00  & 120.00  & 120.00  & 120.00  \\
			12                                      & 120.00  & 120.00  & 120.00  & 120.00  & 120.00  & 120.00  & 120.00  & 120.00  & 120.00  & 120.00  \\
			14                                      & 120.00  & 120.00  & 120.00  & 120.00  & 120.00  & 120.00  & 120.00  & 120.00  & 120.00  & 120.00  \\
			16                                      & 120.00  & 120.00  & 120.00  & 120.00  & 120.00  & 120.00  & 120.00  & 120.00  & 120.00  & 120.00  \\
			\rowcolor{LightCyan}
			18 & 124.44  & 144.44  & 149.66  & 163.12  & 87.03   & 73.35   & 109.48  & 107.28  & 106.19  & 129.00  \\
			\rowcolor{LightCyan}
			19 & 125.14  & 110.91  & 186.74  & 147.19  & 109.15  & 83.72   & 145.52  & 103.08  & 127.51  & 125.47  \\
			\rowcolor{LightCyan}
			20 & 112.67  & 103.97  & 159.03  & 145.64  & 86.44   & 79.05   & 141.10  & 109.85  & 112.92  & 96.27   \\
			\rowcolor{LightCyan}
			21 & 80.00   & 80.00   & 120.00  & 120.00  & 60.00   & 60.00   & 80.00   & 80.00   & 80.00   & 80.00   \\
			\rowcolor{LightCyan}
			22 & 80.00   & 80.00   & 120.00  & 120.00  & 60.00   & 60.00   & 80.00   & 80.00   & 80.00   & 80.00   \\
			\rowcolor{LightCyan}
			23 & 80.00   & 80.00   & 120.00  & 120.00  & 60.00   & 60.00   & 80.00   & 80.00   & 80.00   & 80.00   \\ \hline
		\end{tabular}%
	}
\end{table}

\begin{table}[!tbp]
	\caption{Changes in the average link travel times (in seconds) of AVs in different traffic scenarios (source centroid connectors and detour links are highlighted in yellow and cyan respectively)}
	\label{tab:avttmultipath}
	\resizebox{\textwidth}{!}{%
		\begin{tabular}{lllllllllll}
			\hline
			\multirow{2}{*}{Links} & \multicolumn{2}{c}{Base} & \multicolumn{2}{c}{Longer detour} & \multicolumn{2}{c}{Shorter detour} & \multicolumn{2}{c}{Decreased capacity} & \multicolumn{2}{c}{Increased capacity} \\ \cmidrule(l){2-3} \cmidrule(l){4-5} \cmidrule(l){6-7} \cmidrule(l){8-8} \cmidrule(l){9-9} \cmidrule(l){10-10} \cmidrule(l){11-11} 
			& Initial     & Optimal    & Initial         & Optimal        & Initial         & Optimal        & Initial        & Optimal        & Initial            & Optimal            \\ \midrule
			\rowcolor{yellow}
			1  & 171.67  & 167.39  & 186.91  & 170.60  & 407.15  & 179.76  & 118.66  & 176.48  & 178.19  & 179.43  \\
			2  & 181.17  & 143.64  & 321.62  & 135.51  & 232.02  & 152.23  & 220.47  & 176.77  & 192.38  & 156.04  \\
			\rowcolor{yellow}
			3  & 162.00  & 138.96  & 378.12  & 128.76  & 352.01  & 178.31  & 166.43  & 196.02  & 140.72  & 125.08  \\
			4  & 180.64  & 196.20  & 154.67  & 134.13  & 130.41  & 152.59  & 159.11  & 137.35  & 185.24  & 160.52  \\
			5  & 181.48  & 140.43  & 155.06  & 147.29  & 147.44  & 144.40  & 185.05  & 132.89  & 182.91  & 159.84  \\
			\rowcolor{yellow}
			6  & 307.86  & 116.17  & 159.53  & 123.77  & 672.58  & 151.25  & 128.75  & 155.06  & 210.25  & 152.76  \\
			7  & 127.43  & 122.45  & 153.43  & 127.62  & 127.12  & 126.76  & 154.18  & 124.80  & 154.71  & 128.73  \\
			8  & 128.86  & 133.49  & 136.32  & 128.08  & 129.19  & 128.48  & 130.41  & 129.76  & 127.67  & 125.27  \\
			\rowcolor{yellow}
			10 & 348.59  & 153.93  & 192.02  & 118.59  & 408.32  & 151.30  & 323.37  & 142.01  & 289.36  & 174.22  \\
			11 & 120.00  & 120.53  & 120.00  & 120.29  & 120.00  & 120.27  & 120.00  & 120.44  & 120.00  & 120.58  \\
			12 & 120.00  & 120.00  & 120.00  & 120.00  & 120.00  & 120.00  & 120.00  & 120.00  & 120.00  & 120.00  \\
			14 & 120.00  & 120.00  & 120.00  & 120.00  & 120.00  & 120.00  & 120.00  & 120.00  & 120.00  & 120.00  \\
			16 & 120.00  & 120.00  & 120.00  & 120.00  & 120.00  & 120.00  & 120.00  & 120.00  & 120.00  & 120.00  \\
			\rowcolor{LightCyan}
			18 & 124.05  & 118.32  & 132.99  & 162.68  & 74.26   & 90.34   & 115.61  & 95.27   & 147.74  & 109.12  \\
			\rowcolor{LightCyan}
			19 & 113.04  & 117.08  & 147.21  & 153.82  & 71.59   & 85.40   & 89.65   & 103.47  & 109.52  & 102.55  \\
			\rowcolor{LightCyan}
			20 & 104.95  & 98.16   & 149.57  & 137.41  & 80.21   & 65.96   & 156.13  & 114.42  & 124.98  & 97.79   \\
			\rowcolor{LightCyan}
			21 & 80.00   & 80.00   & 120.00  & 120.00  & 60.00   & 60.00   & 80.00   & 80.00   & 80.00   & 80.00   \\
			\rowcolor{LightCyan}
			22 & 80.00   & 80.00   & 120.00  & 120.00  & 60.00   & 60.00   & 80.00   & 80.00   & 80.00   & 80.00   \\
			\rowcolor{LightCyan}
			23 & 80.00   & 80.00   & 120.00  & 120.00  & 60.00   & 60.00   & 80.00   & 80.00   & 80.00   & 80.00   \\ \hline
		\end{tabular}%
	}
\end{table}
\section{Conclusion}
\label{conclusion}

In this study, we proposed an optimization framework to solve a multi-OD freeway network design problem for optimal lane design under endogenous AV demand. The traffic in the network includes two vehicle classes: LVs without any automated features and AVs with basic automated features such as CACC, speed harmonisation and cooperative merging. Due to the presence of binary lane design variables and the endogenous demand model, the proposed formulation results in a nonconvex MINLP. We tackle this challenging problem by introducing Benders' decomposition approach which iteratively explores possible lane designs in the master problem and at each iteration solves a sequence of SODTA problems which is shown to converge to fixed-points representative of logit-compatible demand splits. The proposed approach is implemented on three hypothetical freeway networks with single and multiple origins and destinations. On the single-OD network, the fixed-point algorithm is found to converge at multiple fixed points providing different proportions of AV demand. However, these multiple fixed-points are found to have no effect on the objective function of the problem. Here, we prove that for a fixed lane design, there exists at least one fixed-point representing the proportion of AV demand in the network. We also prove that the proposed solution method converges to a local optima of the nonconvex problem and identify under which conditions this local optima is a global solution.
{While we cannot guarantee that our algorithm will find the optimal lane design due to the aforementioned non-uniqueness of fixed-points, our numerical experiments have shown that 1) for any lane design there exists at least one fixed-point (see Proposition \ref{prop1}) and 2) that Algorithm \ref{algo:lbfp} converges to a local optima of the FNDP which corresponds to fixed-point solution. Further, we have also provided empirical evidence of the existence of multiple fixed-points (see Section \ref{fp}).}The numerical results on the multi-OD networks show that it is not beneficial in terms of system performance to provide AV lanes for all the links in the network.

\textcolor{blue}{In this study, we adopt the widely used method of successive averages (MSA) to find fixed-point solutions. In our numerical experiments, we were successful in finding fixed point solutions for the proposed formulation, typically within a few MSA iterations. This suggests that although theoretical convergence is not guaranteed, the proposed algorithm is able to find feasible solutions for the Freeway Network Design Problem at hand. Further, as shown in several studies on DTA and SODTA, the MSA algorithm has been widely used to generate solutions to these problems without guaranteed theoretical convergence \cite{wu2003network,bell2008attacker,sbayti2007efficient}}.

\textcolor{blue}
{From a practical application standpoint, the AV-lane allocation problem proposed in this study develops a short term (e.g., morning peak hour) resource allocation strategy and requires minimal infrastructural setup to convert candidate regular lanes to AV-exclusive lanes. This infrastructural setup might include only a variable message sign board fixed at the entrance of each candidate regular lane to inform the road users whether a candidate lane is an AV-exclusive lane or not. The overall system performance metric that involves minimising all vehicles' travel time in the network, decides which candidate regular lane to convert to an AV-exclusive lane during the analysis period. This lane allocation might change during the next analysis period depending on the temporal dynamics of vehicles in the network. Hence, travel time per vehicle type during the analysis period becomes a critical component in allocating the exclusive lanes.
Further, the proposed model can be used for long-term strategic planning models as well. Due to minimal infrastructure investment in setting up AV-exclusive lanes, this approach might allow the freeway to be potentially reconfigured during the off-peak / evening peak. The proposed model may also be useful for designing ramp metering for a freeway network with AV-exclusive lanes.} 

\textcolor{blue}{The following limitations may be examined in future research. We observe that the optimal lane design of freeway network is non-trivial while accounting for endogenous demand of each mode. As shown/discussed in Fig 3/Section 2.3, the proposed formulation may generate lane designs to admit multiple fixed points solutions. Such behaviour is typical of DTA formulations and can be understood via their nonconvex mathematical representations. The presence of the logit-based mode choice model within the proposed formulation for the FNDP corresponds to a nonconvex solution space, even for fixed lane designs. From a practical standpoint, sensitivity analyses around the obtained lane design may be conducted to increase the level of confidence in the solution.} Further, we do not model lane-changing behaviour nor attempt to model vehicle holding issues in this study. {While the proposed model can handle multi-path per OD, the knowledge of a path set is required to determine the average OD travel time. In future research, this limitation could be tackled with a column generation procedure to iteratively augment the path sets.}
This work can also be extended by investigating user equilibrium (UE) vs SO route choice in an arterial network with dedicated AV lanes incorporating mixed vehicular interactions.
\bibliographystyle{elsarticle-harv}
\bibliography{sample}

\begin{thebibliography}{50}
\expandafter\ifx\csname natexlab\endcsname\relax\def\natexlab#1{#1}\fi
\providecommand{\url}[1]{\texttt{#1}}
\providecommand{\href}[2]{#2}
\providecommand{\path}[1]{#1}
\providecommand{\DOIprefix}{doi:}
\providecommand{\ArXivprefix}{arXiv:}
\providecommand{\URLprefix}{URL: }
\providecommand{\Pubmedprefix}{pmid:}
\providecommand{\doi}[1]{\href{http://dx.doi.org/#1}{\path{#1}}}
\providecommand{\Pubmed}[1]{\href{pmid:#1}{\path{#1}}}
\providecommand{\bibinfo}[2]{#2}
\ifx\xfnm\relax \def\xfnm[#1]{\unskip,\space#1}\fi
\bibitem[{Bell et~al.(2008)Bell, Kanturska, Schm{\"o}cker and
  Fonzone}]{bell2008attacker}
\bibinfo{author}{Bell, M.G.}, \bibinfo{author}{Kanturska, U.},
  \bibinfo{author}{Schm{\"o}cker, J.D.}, \bibinfo{author}{Fonzone, A.},
  \bibinfo{year}{2008}.
\newblock \bibinfo{title}{Attacker--defender models and road network
  vulnerability}.
\newblock \bibinfo{journal}{Philosophical Transactions of the Royal Society A:
  Mathematical, Physical and Engineering Sciences} \bibinfo{volume}{366},
  \bibinfo{pages}{1893--1906}.
\bibitem[{Chakraborty et~al.(2018)Chakraborty, Rey, Moylan and
  Waller}]{shantanu2018}
\bibinfo{author}{Chakraborty, S.}, \bibinfo{author}{Rey, D.},
  \bibinfo{author}{Moylan, E.}, \bibinfo{author}{Waller, S.T.},
  \bibinfo{year}{2018}.
\newblock \bibinfo{title}{Link transmission model-based linear programming
  formulation for network design}.
\newblock \bibinfo{journal}{Transportation Research Record}
  \bibinfo{volume}{2672}, \bibinfo{pages}{139--147}.
\bibitem[{Chen et~al.(2016)Chen, He, Zhang and Yin}]{chen2016optimal}
\bibinfo{author}{Chen, Z.}, \bibinfo{author}{He, F.}, \bibinfo{author}{Zhang,
  L.}, \bibinfo{author}{Yin, Y.}, \bibinfo{year}{2016}.
\newblock \bibinfo{title}{Optimal deployment of autonomous vehicle lanes with
  endogenous market penetration}.
\newblock \bibinfo{journal}{Transportation Research Part C: Emerging
  Technologies} \bibinfo{volume}{72}, \bibinfo{pages}{143--156}.
\bibitem[{Courant et~al.(1928)Courant, Friedrichs and
  Lewy}]{courant1928partiellen}
\bibinfo{author}{Courant, R.}, \bibinfo{author}{Friedrichs, K.},
  \bibinfo{author}{Lewy, H.}, \bibinfo{year}{1928}.
\newblock \bibinfo{title}{{\"U}ber die partiellen differenzengleichungen der
  mathematischen physik}.
\newblock \bibinfo{journal}{Mathematische annalen} \bibinfo{volume}{100},
  \bibinfo{pages}{32--74}.
\bibitem[{Daganzo(1994)}]{daganzo1994cell}
\bibinfo{author}{Daganzo, C.F.}, \bibinfo{year}{1994}.
\newblock \bibinfo{title}{The cell transmission model: A dynamic representation
  of highway traffic consistent with the hydrodynamic theory}.
\newblock \bibinfo{journal}{Transportation Research Part B: Methodological}
  \bibinfo{volume}{28}, \bibinfo{pages}{269--287}.
\bibitem[{Dresner and Stone(2004)}]{dresner2004multiagent}
\bibinfo{author}{Dresner, K.}, \bibinfo{author}{Stone, P.},
  \bibinfo{year}{2004}.
\newblock \bibinfo{title}{Multiagent traffic management: A reservation-based
  intersection control mechanism}, in: \bibinfo{booktitle}{Proceedings of the
  Third International Joint Conference on Autonomous Agents and Multiagent
  Systems-Volume 2}, \bibinfo{organization}{IEEE Computer Society}. pp.
  \bibinfo{pages}{530--537}.
\bibitem[{Fajardo et~al.(2011)Fajardo, Au, Waller, Stone and
  Yang}]{fajardo2011automated}
\bibinfo{author}{Fajardo, D.}, \bibinfo{author}{Au, T.C.},
  \bibinfo{author}{Waller, S.T.}, \bibinfo{author}{Stone, P.},
  \bibinfo{author}{Yang, D.}, \bibinfo{year}{2011}.
\newblock \bibinfo{title}{Automated intersection control: Performance of future
  innovation versus current traffic signal control}.
\newblock \bibinfo{journal}{Transportation Research Record}
  \bibinfo{volume}{2259}, \bibinfo{pages}{223--232}.
\bibitem[{Fernandes and Nunes(2012)}]{fernandes2012platooning}
\bibinfo{author}{Fernandes, P.}, \bibinfo{author}{Nunes, U.},
  \bibinfo{year}{2012}.
\newblock \bibinfo{title}{Platooning with ivc-enabled autonomous vehicles:
  Strategies to mitigate communication delays, improve safety and traffic
  flow}.
\newblock \bibinfo{journal}{IEEE Transactions on Intelligent Transportation
  Systems} \bibinfo{volume}{13}, \bibinfo{pages}{91--106}.
\bibitem[{Greenblatt and Saxena(2015)}]{greenblatt2015autonomous}
\bibinfo{author}{Greenblatt, J.B.}, \bibinfo{author}{Saxena, S.},
  \bibinfo{year}{2015}.
\newblock \bibinfo{title}{Autonomous taxis could greatly reduce greenhouse-gas
  emissions of us light-duty vehicles}.
\newblock \bibinfo{journal}{Nature Climate Change} \bibinfo{volume}{5},
  \bibinfo{pages}{860}.
\bibitem[{Greenshields et~al.(1935)Greenshields, Bibbins, Channing and
  Miller}]{greenshields}
\bibinfo{author}{Greenshields, B.}, \bibinfo{author}{Bibbins, J.},
  \bibinfo{author}{Channing, W.}, \bibinfo{author}{Miller, H.},
  \bibinfo{year}{1935}.
\newblock \bibinfo{title}{A study of traffic capacity}.
\newblock \bibinfo{journal}{Highway Research Board proceedings} .
\bibitem[{Haboucha et~al.(2017)Haboucha, Ishaq and Shiftan}]{haboucha2017user}
\bibinfo{author}{Haboucha, C.J.}, \bibinfo{author}{Ishaq, R.},
  \bibinfo{author}{Shiftan, Y.}, \bibinfo{year}{2017}.
\newblock \bibinfo{title}{User preferences regarding autonomous vehicles}.
\newblock \bibinfo{journal}{Transportation Research Part C: Emerging
  Technologies} \bibinfo{volume}{78}, \bibinfo{pages}{37--49}.
\bibitem[{Hyland and Mahmassani(2018)}]{hyland2018dynamic}
\bibinfo{author}{Hyland, M.}, \bibinfo{author}{Mahmassani, H.S.},
  \bibinfo{year}{2018}.
\newblock \bibinfo{title}{Dynamic autonomous vehicle fleet operations:
  Optimization-based strategies to assign avs to immediate traveler demand
  requests}.
\newblock \bibinfo{journal}{Transportation Research Part C: Emerging
  Technologies} \bibinfo{volume}{92}, \bibinfo{pages}{278--297}.
\bibitem[{Krueger et~al.(2016)Krueger, Rashidi and
  Rose}]{krueger2016preferences}
\bibinfo{author}{Krueger, R.}, \bibinfo{author}{Rashidi, T.H.},
  \bibinfo{author}{Rose, J.M.}, \bibinfo{year}{2016}.
\newblock \bibinfo{title}{Preferences for shared autonomous vehicles}.
\newblock \bibinfo{journal}{Transportation research part C: emerging
  technologies} \bibinfo{volume}{69}, \bibinfo{pages}{343--355}.
\bibitem[{Levin(2017)}]{levin2017congestion}
\bibinfo{author}{Levin, M.W.}, \bibinfo{year}{2017}.
\newblock \bibinfo{title}{Congestion-aware system optimal route choice for
  shared autonomous vehicles}.
\newblock \bibinfo{journal}{Transportation Research Part C: Emerging
  Technologies} \bibinfo{volume}{82}, \bibinfo{pages}{229--247}.
\bibitem[{Levin and Boyles(2016a)}]{levin2016cell}
\bibinfo{author}{Levin, M.W.}, \bibinfo{author}{Boyles, S.D.},
  \bibinfo{year}{2016}a.
\newblock \bibinfo{title}{A cell transmission model for dynamic lane reversal
  with autonomous vehicles}.
\newblock \bibinfo{journal}{Transportation Research Part C: Emerging
  Technologies} \bibinfo{volume}{68}, \bibinfo{pages}{126--143}.
\bibitem[{Levin and Boyles(2016b)}]{LEVIN2016103}
\bibinfo{author}{Levin, M.W.}, \bibinfo{author}{Boyles, S.D.},
  \bibinfo{year}{2016}b.
\newblock \bibinfo{title}{A multiclass cell transmission model for shared human
  and autonomous vehicle roads}.
\newblock \bibinfo{journal}{Transportation Research Part C: Emerging
  Technologies} \bibinfo{volume}{62}, \bibinfo{pages}{103 -- 116}.
\bibitem[{Levin and Rey(2017)}]{levin2017conflict}
\bibinfo{author}{Levin, M.W.}, \bibinfo{author}{Rey, D.}, \bibinfo{year}{2017}.
\newblock \bibinfo{title}{Conflict-point formulation of intersection control
  for autonomous vehicles}.
\newblock \bibinfo{journal}{Transportation Research Part C: Emerging
  Technologies} \bibinfo{volume}{85}, \bibinfo{pages}{528--547}.
\bibitem[{Lighthill and Whitham(1955)}]{lighthill1955kinematic}
\bibinfo{author}{Lighthill, M.J.}, \bibinfo{author}{Whitham, G.B.},
  \bibinfo{year}{1955}.
\newblock \bibinfo{title}{On kinematic waves ii. a theory of traffic flow on
  long crowded roads}.
\newblock \bibinfo{journal}{Proceedings of the Royal Society of London. Series
  A. Mathematical and Physical Sciences} \bibinfo{volume}{229},
  \bibinfo{pages}{317--345}.
\bibitem[{Long et~al.(2018)Long, Chen, Szeto and Shi}]{long2018link}
\bibinfo{author}{Long, J.}, \bibinfo{author}{Chen, J.}, \bibinfo{author}{Szeto,
  W.}, \bibinfo{author}{Shi, Q.}, \bibinfo{year}{2018}.
\newblock \bibinfo{title}{Link-based system optimum dynamic traffic assignment
  problems with environmental objectives}.
\newblock \bibinfo{journal}{Transportation Research Part D: Transport and
  Environment} \bibinfo{volume}{60}, \bibinfo{pages}{56--75}.
\bibitem[{Long and Szeto(2019)}]{long2019link}
\bibinfo{author}{Long, J.}, \bibinfo{author}{Szeto, W.Y.},
  \bibinfo{year}{2019}.
\newblock \bibinfo{title}{Link-based system optimum dynamic traffic assignment
  problems in general networks}.
\newblock \bibinfo{journal}{Operations Research} \bibinfo{volume}{67},
  \bibinfo{pages}{167--182}.
\bibitem[{Ma and Wang(2019)}]{ma2019influence}
\bibinfo{author}{Ma, K.}, \bibinfo{author}{Wang, H.}, \bibinfo{year}{2019}.
\newblock \bibinfo{title}{Influence of exclusive lanes for connected and
  autonomous vehicles on freeway traffic flow}.
\newblock \bibinfo{journal}{IEEE Access} \bibinfo{volume}{7},
  \bibinfo{pages}{50168--50178}.
\bibitem[{Mahmassani(2016)}]{mahmassani201650th}
\bibinfo{author}{Mahmassani, H.S.}, \bibinfo{year}{2016}.
\newblock \bibinfo{title}{50th anniversary invited article—autonomous
  vehicles and connected vehicle systems: flow and operations considerations}.
\newblock \bibinfo{journal}{Transportation Science} \bibinfo{volume}{50},
  \bibinfo{pages}{1140--1162}.
\bibitem[{Melson et~al.(2018)Melson, Levin, Hammit and Boyles}]{MELSON2018114}
\bibinfo{author}{Melson, C.L.}, \bibinfo{author}{Levin, M.W.},
  \bibinfo{author}{Hammit, B.E.}, \bibinfo{author}{Boyles, S.D.},
  \bibinfo{year}{2018}.
\newblock \bibinfo{title}{Dynamic traffic assignment of cooperative adaptive
  cruise control}.
\newblock \bibinfo{journal}{Transportation Research Part C: Emerging
  Technologies} \bibinfo{volume}{90}, \bibinfo{pages}{114 -- 133}.
\bibitem[{Mersky and Samaras(2016)}]{mersky2016fuel}
\bibinfo{author}{Mersky, A.C.}, \bibinfo{author}{Samaras, C.},
  \bibinfo{year}{2016}.
\newblock \bibinfo{title}{Fuel economy testing of autonomous vehicles}.
\newblock \bibinfo{journal}{Transportation Research Part C: Emerging
  Technologies} \bibinfo{volume}{65}, \bibinfo{pages}{31--48}.
\bibitem[{Milan{\'e}s et~al.(2013)Milan{\'e}s, Shladover, Spring, Nowakowski,
  Kawazoe and Nakamura}]{milanes2013cooperative}
\bibinfo{author}{Milan{\'e}s, V.}, \bibinfo{author}{Shladover, S.E.},
  \bibinfo{author}{Spring, J.}, \bibinfo{author}{Nowakowski, C.},
  \bibinfo{author}{Kawazoe, H.}, \bibinfo{author}{Nakamura, M.},
  \bibinfo{year}{2013}.
\newblock \bibinfo{title}{Cooperative adaptive cruise control in real traffic
  situations}.
\newblock \bibinfo{journal}{IEEE Transactions on Intelligent Transportation
  Systems} \bibinfo{volume}{15}, \bibinfo{pages}{296--305}.
\bibitem[{Monteil et~al.(2018)Monteil, Bouroche and Leith}]{monteil2018mathcal}
\bibinfo{author}{Monteil, J.}, \bibinfo{author}{Bouroche, M.},
  \bibinfo{author}{Leith, D.J.}, \bibinfo{year}{2018}.
\newblock \bibinfo{title}{Stability analysis of heterogeneous traffic with
  application to parameter optimization for the control of automated vehicles}.
\newblock \bibinfo{journal}{IEEE Transactions on Control Systems Technology}
  \bibinfo{volume}{27}, \bibinfo{pages}{934--949}.
\bibitem[{Movaghar et~al.(2020)Movaghar, Mesbah and
  Habibian}]{movaghar2020optimum}
\bibinfo{author}{Movaghar, S.}, \bibinfo{author}{Mesbah, M.},
  \bibinfo{author}{Habibian, M.}, \bibinfo{year}{2020}.
\newblock \bibinfo{title}{Optimum location of autonomous vehicle lanes: A model
  considering capacity variation}.
\newblock \bibinfo{journal}{Mathematical Problems in Engineering}
  \bibinfo{volume}{2020}.
\bibitem[{Newell(1993)}]{newell1993simplified}
\bibinfo{author}{Newell, G.F.}, \bibinfo{year}{1993}.
\newblock \bibinfo{title}{A simplified theory of kinematic waves in highway
  traffic, part ii: Queueing at freeway bottlenecks}.
\newblock \bibinfo{journal}{Transportation Research Part B: Methodological}
  \bibinfo{volume}{27}, \bibinfo{pages}{289--303}.
\bibitem[{Ploeg et~al.(2011)Ploeg, Scheepers, Van~Nunen, Van~de Wouw and
  Nijmeijer}]{ploeg2011design}
\bibinfo{author}{Ploeg, J.}, \bibinfo{author}{Scheepers, B.T.},
  \bibinfo{author}{Van~Nunen, E.}, \bibinfo{author}{Van~de Wouw, N.},
  \bibinfo{author}{Nijmeijer, H.}, \bibinfo{year}{2011}.
\newblock \bibinfo{title}{Design and experimental evaluation of cooperative
  adaptive cruise control}, in: \bibinfo{booktitle}{2011 14th International
  IEEE Conference on Intelligent Transportation Systems (ITSC)},
  \bibinfo{organization}{IEEE}. pp. \bibinfo{pages}{260--265}.
\bibitem[{Powell and Sheffi(1982)}]{powell1982convergence}
\bibinfo{author}{Powell, W.B.}, \bibinfo{author}{Sheffi, Y.},
  \bibinfo{year}{1982}.
\newblock \bibinfo{title}{The convergence of equilibrium algorithms with
  predetermined step sizes}.
\newblock \bibinfo{journal}{Transportation Science} \bibinfo{volume}{16},
  \bibinfo{pages}{45--55}.
\bibitem[{Qian et~al.(2014)Qian, Gregoire, Moutarde and
  De~La~Fortelle}]{qian2014priority}
\bibinfo{author}{Qian, X.}, \bibinfo{author}{Gregoire, J.},
  \bibinfo{author}{Moutarde, F.}, \bibinfo{author}{De~La~Fortelle, A.},
  \bibinfo{year}{2014}.
\newblock \bibinfo{title}{Priority-based coordination of autonomous and legacy
  vehicles at intersection}, in: \bibinfo{booktitle}{17th international IEEE
  conference on intelligent transportation systems (ITSC)},
  \bibinfo{organization}{IEEE}. pp. \bibinfo{pages}{1166--1171}.
\bibitem[{Rahmaniani et~al.(2017)Rahmaniani, Crainic, Gendreau and
  Rei}]{RAHMANIANI2017801}
\bibinfo{author}{Rahmaniani, R.}, \bibinfo{author}{Crainic, T.G.},
  \bibinfo{author}{Gendreau, M.}, \bibinfo{author}{Rei, W.},
  \bibinfo{year}{2017}.
\newblock \bibinfo{title}{The benders decomposition algorithm: A literature
  review}.
\newblock \bibinfo{journal}{European Journal of Operational Research}
  \bibinfo{volume}{259}, \bibinfo{pages}{801 -- 817}.
\bibitem[{Rey and Levin(2019)}]{rey2019blue}
\bibinfo{author}{Rey, D.}, \bibinfo{author}{Levin, M.W.}, \bibinfo{year}{2019}.
\newblock \bibinfo{title}{Blue phase: Optimal network traffic control for
  legacy and autonomous vehicles}.
\newblock \bibinfo{journal}{Transportation Research Part B: Methodological}
  \bibinfo{volume}{130}, \bibinfo{pages}{105--129}.
\bibitem[{Richards(1956)}]{richards1956shock}
\bibinfo{author}{Richards, P.I.}, \bibinfo{year}{1956}.
\newblock \bibinfo{title}{Shock waves on the highway}.
\newblock \bibinfo{journal}{Operations research} \bibinfo{volume}{4},
  \bibinfo{pages}{42--51}.
\bibitem[{SAE(2013)}]{sae2013definitions}
\bibinfo{author}{SAE, T.}, \bibinfo{year}{2013}.
\newblock \bibinfo{title}{Definitions for terms related to on-road motor
  vehicle automated driving systems-j3016}.
\newblock \bibinfo{journal}{Society of Automotive Engineers: On-Road Automated
  Vehicle Standards Committee; SAE Pub. Inc., Warrendale, PA, USA} .
\bibitem[{Sbayti et~al.(2007)Sbayti, Lu and Mahmassani}]{sbayti2007efficient}
\bibinfo{author}{Sbayti, H.}, \bibinfo{author}{Lu, C.C.},
  \bibinfo{author}{Mahmassani, H.S.}, \bibinfo{year}{2007}.
\newblock \bibinfo{title}{Efficient implementation of method of successive
  averages in simulation-based dynamic traffic assignment models for
  large-scale network applications}.
\newblock \bibinfo{journal}{Transportation Research Record}
  \bibinfo{volume}{2029}, \bibinfo{pages}{22--30}.
\bibitem[{Sheffi(1985)}]{sheffi1985urban}
\bibinfo{author}{Sheffi, Y.}, \bibinfo{year}{1985}.
\newblock \bibinfo{title}{Urban transportation networks}.
  volume~\bibinfo{volume}{6}.
\newblock \bibinfo{publisher}{Prentice-Hall, Englewood Cliffs, NJ}.
\bibitem[{Shladover et~al.(2012)Shladover, Su and Lu}]{shladover2012impacts}
\bibinfo{author}{Shladover, S.E.}, \bibinfo{author}{Su, D.},
  \bibinfo{author}{Lu, X.Y.}, \bibinfo{year}{2012}.
\newblock \bibinfo{title}{Impacts of cooperative adaptive cruise control on
  freeway traffic flow}.
\newblock \bibinfo{journal}{Transportation Research Record}
  \bibinfo{volume}{2324}, \bibinfo{pages}{63--70}.
\bibitem[{Sullivan(2018)}]{frost&sullivan}
\bibinfo{author}{Sullivan, F..}, \bibinfo{year}{2018}.
\newblock \bibinfo{title}{Global autonomous driving market outlook}.
\newblock \bibinfo{journal}{Online Store} \bibinfo{volume}{K24A}.
\bibitem[{Talebpour and Mahmassani(2016)}]{talebpour2016influence}
\bibinfo{author}{Talebpour, A.}, \bibinfo{author}{Mahmassani, H.S.},
  \bibinfo{year}{2016}.
\newblock \bibinfo{title}{Influence of connected and autonomous vehicles on
  traffic flow stability and throughput}.
\newblock \bibinfo{journal}{Transportation Research Part C: Emerging
  Technologies} \bibinfo{volume}{71}, \bibinfo{pages}{143--163}.
\bibitem[{Talebpour et~al.(2017)Talebpour, Mahmassani and Elfar}]{alireza2017}
\bibinfo{author}{Talebpour, A.}, \bibinfo{author}{Mahmassani, H.S.},
  \bibinfo{author}{Elfar, A.}, \bibinfo{year}{2017}.
\newblock \bibinfo{title}{Investigating the effects of reserved lanes for
  autonomous vehicles on congestion and travel time reliability}.
\newblock \bibinfo{journal}{Transportation Research Record}
  \bibinfo{volume}{2622}, \bibinfo{pages}{1--12}.
\bibitem[{Tientrakool et~al.(2011)Tientrakool, Ho and
  Maxemchuk}]{tientrakool2011highway}
\bibinfo{author}{Tientrakool, P.}, \bibinfo{author}{Ho, Y.C.},
  \bibinfo{author}{Maxemchuk, N.F.}, \bibinfo{year}{2011}.
\newblock \bibinfo{title}{Highway capacity benefits from using
  vehicle-to-vehicle communication and sensors for collision avoidance}, in:
  \bibinfo{booktitle}{2011 IEEE Vehicular Technology Conference (VTC Fall)},
  \bibinfo{organization}{IEEE}. pp. \bibinfo{pages}{1--5}.
\bibitem[{Van~Arem et~al.(2006)Van~Arem, Van~Driel and Visser}]{van2006impact}
\bibinfo{author}{Van~Arem, B.}, \bibinfo{author}{Van~Driel, C.J.},
  \bibinfo{author}{Visser, R.}, \bibinfo{year}{2006}.
\newblock \bibinfo{title}{The impact of cooperative adaptive cruise control on
  traffic-flow characteristics}.
\newblock \bibinfo{journal}{IEEE Transactions on intelligent transportation
  systems} \bibinfo{volume}{7}, \bibinfo{pages}{429--436}.
\bibitem[{Vander~Laan and Sadabadi(2017)}]{vander2017operational}
\bibinfo{author}{Vander~Laan, Z.}, \bibinfo{author}{Sadabadi, K.F.},
  \bibinfo{year}{2017}.
\newblock \bibinfo{title}{Operational performance of a congested corridor with
  lanes dedicated to autonomous vehicle traffic}.
\newblock \bibinfo{journal}{International Journal of Transportation Science and
  Technology} \bibinfo{volume}{6}, \bibinfo{pages}{42--52}.
\bibitem[{Wang et~al.(2019)Wang, Peeta and He}]{wang2019multiclass}
\bibinfo{author}{Wang, J.}, \bibinfo{author}{Peeta, S.}, \bibinfo{author}{He,
  X.}, \bibinfo{year}{2019}.
\newblock \bibinfo{title}{Multiclass traffic assignment model for mixed traffic
  flow of human-driven vehicles and connected and autonomous vehicles}.
\newblock \bibinfo{journal}{Transportation Research Part B: Methodological}
  \bibinfo{volume}{126}, \bibinfo{pages}{139--168}.
\bibitem[{Wong et~al.(2018)Wong, Saxena and Dixit}]{neeraj2018}
\bibinfo{author}{Wong, T.W.}, \bibinfo{author}{Saxena, N.},
  \bibinfo{author}{Dixit, V.V.}, \bibinfo{year}{2018}.
\newblock \bibinfo{title}{A study of route choice behavior of drivers in
  autonomous vehicles}, in: \bibinfo{booktitle}{Transportation Research Board
  97th Annual Meeting}, \bibinfo{organization}{Transportation Research Board}.
\bibitem[{Wu and Lam(2003)}]{wu2003network}
\bibinfo{author}{Wu, Z.}, \bibinfo{author}{Lam, W.H.}, \bibinfo{year}{2003}.
\newblock \bibinfo{title}{Network equilibrium for congested multi-mode networks
  with elastic demand}.
\newblock \bibinfo{journal}{Journal of Advanced Transportation}
  \bibinfo{volume}{37}, \bibinfo{pages}{295--318}.
\bibitem[{Ye and Yamamoto(2018)}]{ye2018impact}
\bibinfo{author}{Ye, L.}, \bibinfo{author}{Yamamoto, T.}, \bibinfo{year}{2018}.
\newblock \bibinfo{title}{Impact of dedicated lanes for connected and
  autonomous vehicle on traffic flow throughput}.
\newblock \bibinfo{journal}{Physica A: Statistical Mechanics and its
  Applications} \bibinfo{volume}{512}, \bibinfo{pages}{588--597}.
\bibitem[{Yperman et~al.(2005)Yperman, Logghe and Immers}]{yperman2005link}
\bibinfo{author}{Yperman, I.}, \bibinfo{author}{Logghe, S.},
  \bibinfo{author}{Immers, B.}, \bibinfo{year}{2005}.
\newblock \bibinfo{title}{The link transmission model: An efficient
  implementation of the kinematic wave theory in traffic networks}, in:
  \bibinfo{booktitle}{Proceedings of the 10th EWGT Meeting},
  \bibinfo{organization}{Poznan Poland}. pp. \bibinfo{pages}{122--127}.
\bibitem[{Yu et~al.(2019)Yu, Tak, Park and Yeo}]{yu2019impact}
\bibinfo{author}{Yu, H.}, \bibinfo{author}{Tak, S.}, \bibinfo{author}{Park,
  M.}, \bibinfo{author}{Yeo, H.}, \bibinfo{year}{2019}.
\newblock \bibinfo{title}{Impact of autonomous-vehicle-only lanes in mixed
  traffic conditions}.
\newblock \bibinfo{journal}{Transportation Research Record} ,
  \bibinfo{pages}{0361198119847475}.

\end{thebibliography}

\end{document}